\documentclass[12pt]{amsart}
\begin{document}
\theoremstyle{plain}
\newtheorem{thm}{Theorem}[section]
\newtheorem*{thm*}{Theorem}
\newtheorem{prop}[thm]{Proposition}
\newtheorem*{prop*}{Proposition}
\newtheorem{lemma}[thm]{Lemma}
\newtheorem{cor}[thm]{Corollary}
\newtheorem*{conj*}{Conjecture}
\newtheorem*{cor*}{Corollary}
\newtheorem{defn}[thm]{Definition}
\theoremstyle{definition}
\newtheorem*{defn*}{Definition}
\newtheorem{rems}[thm]{Remarks}
\newtheorem*{rems*}{Remarks}
\newtheorem{rem}[thm]{Remark}
\newtheorem*{rem*}{Remark}
\newtheorem*{proof*}{Proof}
\newtheorem*{not*}{Notation}
\newcommand{\npartial}{\slash\!\!\!\partial}
\newcommand{\Heis}{\operatorname{Heis}}
\newcommand{\Solv}{\operatorname{Solv}}
\newcommand{\Spin}{\operatorname{Spin}}
\newcommand{\SO}{\operatorname{SO}}
\newcommand{\ind}{\operatorname{ind}}
\newcommand{\Index}{\operatorname{index}}
\newcommand{\ch}{\operatorname{ch}}
\newcommand{\rank}{\operatorname{rank}}
\newcommand{\abs}[1]{\lvert#1\rvert}
 \newcommand{\A}{{\mathcal A}}
        \newcommand{\D}{{\mathcal D}}\newcommand{\HH}{{\mathcal H}}
        \newcommand{\LL}{{\mathcal L}}
        \newcommand{\B}{{\mathcal B}}
        \newcommand{\K}{{\mathcal K}}
\newcommand{\oo}{{\mathcal O}}
         \newcommand{\PP}{{\mathcal P}}
        \newcommand{\s}{\sigma}
        \newcommand{\coker}{{\mbox coker}}
        \newcommand{\p}{\partial}
        \newcommand{\dd}{|\D|}
        \newcommand{\n}{\Vert}  
\newcommand{\bma}{\left(\begin{array}{cc}}
\newcommand{\ema}{\end{array}\right)}
\newcommand{\bca}{\left(\begin{array}{c}}
\newcommand{\eca}{\end{array}\right)}

\newcommand{\sr}{\stackrel}
\newcommand{\da}{\downarrow}
\newcommand{\tD}{\tilde{\D}}

        \newcommand{\R}{\mathbf R}
        \newcommand{\C}{\mathbf C}
        \newcommand{\h}{\mathbf H}
\newcommand{\Z}{\mathbf Z}
\newcommand{\N}{\mathbf N}
\newcommand{\tto}{\longrightarrow}
\newcommand{\ben}{\begin{displaymath}}
        \newcommand{\een}{\end{displaymath}}
\newcommand{\be}{\begin{equation}}
\newcommand{\ee}{\end{equation}}

        \newcommand{\bean}{\begin{eqnarray*}}
        \newcommand{\eean}{\end{eqnarray*}}
\newcommand{\nno}{\nonumber\\}
\newcommand{\bea}{\begin{eqnarray}}
        \newcommand{\eea}{\end{eqnarray}}












\newcommand{\supp}[1]{\operatorname{#1}}
\newcommand{\norm}[1]{\parallel\, #1\, \parallel}
\newcommand{\ip}[2]{\langle #1,#2\rangle}
\setlength{\parskip}{.3cm}
\newcommand{\nc}{\newcommand}
\nc{\nt}{\newtheorem}
\nc{\gf}[2]{\genfrac{}{}{0pt}{}{#1}{#2}}
\nc{\mb}[1]{{\mbox{$ #1 $}}}
\nc{\real}{{\mathbb R}}
\nc{\comp}{{\mathbb C}}
\nc{\ints}{{\mathbb Z}}
\nc{\Ltoo}{\mb{L^2({\mathbf H})}}
\nc{\rtoo}{\mb{{\mathbf R}^2}}
\nc{\slr}{{\mathbf {SL}}(2,\real)}
\nc{\slz}{{\mathbf {SL}}(2,\ints)}
\nc{\su}{{\mathbf {SU}}(1,1)}
\nc{\so}{{\mathbf {SO}}}
\nc{\hyp}{{\mathbb H}}
\nc{\disc}{{\mathbf D}}
\nc{\torus}{{\mathbb T}}
\newcommand{\tk}{\widetilde{K}}
\newcommand{\boe}{{\bf e}}\newcommand{\bt}{{\bf t}}
\newcommand{\vth}{\vartheta}
\newcommand{\CGh}{\widetilde{\CG}}
\newcommand{\db}{\overline{\partial}}
\newcommand{\tE}{\widetilde{E}}
\newcommand{\tr}{\mbox{tr}}
\newcommand{\ta}{\widetilde{\alpha}}
\newcommand{\tb}{\widetilde{\beta}}
\newcommand{\txi}{\widetilde{\xi}}
\newcommand{\hV}{\hat{V}}
\newcommand{\IC}{\mathbf{C}}
\newcommand{\IZ}{\mathbf{Z}}
\newcommand{\IP}{\mathbf{P}}
\newcommand{\IR}{\mathbf{R}}
\newcommand{\IH}{\mathbf{H}}
\newcommand{\IG}{\mathbf{G}}
\newcommand{\CC}{{\mathcal C}}
\newcommand{\CD}{{\mathcal D}}
\newcommand{\CS}{{\mathcal S}}
\newcommand{\CG}{{\mathcal G}}
\newcommand{\CL}{{\mathcal L}}
\newcommand{\CO}{{\mathcal O}}
\nc{\ca}{{\mathcal A}}
\nc{\cag}{{{\mathcal A}^\Gamma}}
\nc{\cg}{{\mathcal G}}
\nc{\chh}{{\mathcal H}}
\nc{\ck}{{\mathcal B}}
\nc{\cl}{{\mathcal L}}
\nc{\cm}{{\mathcal M}}
\nc{\cn}{{\mathcal N}}
\nc{\NN}{{\mathcal N}}
\nc{\cs}{{\mathcal S}}
\nc{\cz}{{\mathcal Z}}
\nc{\br}{{\mathcal R}}
\nc{\sind}{\sigma{\rm -ind}}
\newcommand{\la}{\langle}
\newcommand{\ra}{\rangle}

\begin{center}
 \title{The  Chern Character of Semifinite
Spectral Triples}

 \vspace{.5 in}

\author{}

\maketitle

{\bf Alan L. Carey}\\
Mathematical Sciences Institute\\
Australian National University\\
Canberra, ACT. 0200, AUSTRALIA\\
e-mail: acarey@maths.anu.edu.au\\ 
\vspace{.2 in}

{\bf John Phillips}\\
Department of Mathematics and Statistics\\
University of Victoria\\Victoria, B.C. V8W 3P4, CANADA\footnote{Address for
correspondence}\\
e-mail: phillips@math.uvic.ca\\ 
\vspace{.2 in}

{\bf Adam Rennie}\\
Institute for Mathematical Sciences\\University of Copenhagen\\
Universitetsparken 5, DK-2100, Copenhagen, DENMARK\\
e-mail: rennie@math.ku.dk

\vspace{.2 in}

{\bf Fyodor A. Sukochev}\\
School of Informatics and Engineering\\
Flinders University\\
Bedford Park S.A 5042 AUSTRALIA\\
e-mail: sukochev@infoeng.flinders.edu.au\\ 

\vspace{.25 in}

All authors were supported by grants from ARC (Australia) and
NSERC (Canada), in addition the third named author acknowledges a 
University of Newcastle early career researcher grant, and support
from the SNF, Denmark.

\end{center}

\newpage
\centerline{{\bf Abstract}}
In previous work we
 generalised both the odd and 
even local index formula of Connes and Moscovici to the
case of spectral triples for a $*$-subalgebra $\A$
of a general semifinite von Neumann algebra. 
Our proofs are novel even in the setting of the original theorem 
and rely on the introduction of a function valued cocycle 
(called the resolvent cocycle)
which is `almost'
a $(b,B)$-cocycle in the cyclic cohomology of $\mathcal A$. 
In this paper we show that this 
resolvent cocycle `almost' represents the Chern character, and assuming 
analytic continuation properties for zeta functions, we show that the 
associated residue cocycle, which appears in our statement of the 
local index theorem does represent the Chern character. 
\footnote{AMS Subject classification:
Primary: 19K56, 46L80; secondary: 58B30, 46L87. Keywords and Phrases:
von Neumann algebra, Fredholm module, cyclic cohomology, chern character,
spectral flow.}

\newpage
\allowdisplaybreaks
\section{Introduction}

The local index theorem in noncommutative geometry is due to Connes and 
Moscovici \cite{CM}. The main consequence of the theorem is a formula
for the Chern character of an unbounded Fredholm module (or spectral triple)
in terms of  a `residue cocycle' in the $(b,B)$ bicomplex. 
This residue cocycle is a sum of residues of certain zeta functions.
There have been
two new proofs of this formula discovered recently, by Higson 
\cite{hig} and by the present authors \cite{CPRS2, CPRS3}.
The main new feature of these proofs is that the starting point for 
the Connes-Moscovici argument (the JLO cocycle \cite{JLO}, \cite{Co1})
is replaced by a different cocycle derived from  the resolvent expansion
in perturbation theory.
These new proofs have some conceptual and technical advantages over the 
earlier proof. In particular \cite{CPRS2, CPRS3} enable 
the local index theorem to be
extended to the case where one has a spectral triple `inside' a general 
semifinite von Neumann algebra. This extension enabled us to encompass 
examples of differential operators that are not Fredholm in the ordinary 
sense but are Breuer-Fredholm, some examples being described in 
\cite{BeF,CP1,Sh,L,M,BCPRSW}.

On of the  novel features of our proof 
of the semifinite version of the local index formula
is the introduction of a new cocycle
(in the $(b,B)$ bicomplex)
which provides a substitute for the JLO cocycle \cite{Co1}, \cite{JLO} for
finitely summable spectral triples. This 
new cocycle is holomorphic function
valued. 
Similarly Higson introduces an
`improper cocycle' \cite{hig} which is function valued but has 
very different holomorphy properties to our resolvent cocycle.
In the last Section of this paper 
we clarify the relationship of our approach to 
that of Higson. 

What is currently absent from earlier work is a detailed
exposition of the properties of our resolvent cocycle which would 
put it on the same footing as the JLO coycle.
The object of this paper is to 
remedy this situation.
To that end we provide a uniform exposition of both the odd and even
forms of the cocycle that was introduced separately
in each of \cite{CPRS2} (for the odd case) and \cite {CPRS3} 
(for the even case).
Then we give a proof that our resolvent cocycle represents 
(in a suitable sense) 
the Chern character in cyclic cohomology. From this we derive
some important consequences which we will summarise later in the
introduction.

The hypotheses for our proof of the local index formula
are weaker than those of \cite{CM}. 
We assume that we have a finitely summable spectral triple 
$(\ca, \chh, \CD)$ with spectral dimension $q$
(that is $(1+\D^2)^{-n/2}$ is trace class for all $n>q$ and $q$ is the least 
positive real number for which this is true).
This suffices to verify that the
individual functionals that make up the resolvent cocycle
are continuous in an appropriate sense 
with values in functions 
defined and holomorphic in a certain half-plane. We 
obtain a cocycle in the finite 
$(b,B)$ complex by considering the resolvent cocycle modulo those functions 
holomorphic in a half-plane containing the critical point $r=(1-q)/2$
(where {\it a priori} the terms in the resolvent cocycle may have a singularity
and where we take residues to obtain our version of the
 Connes-Moscovici residue cocycle).
We remark, but do not prove here, that the resolvent cocycle is not entire
but at no point do we need infinitely many terms in our
expression for this cocycle.
We use the notation $(\phi^r_m)$, $m=0,1,2,\cdots$ to denote the 
components of the
resolvent cocycle in the $(b,B)$ complex ($m$ odd in the odd case, even in the 
even case).

In our proof of the local index formula in \cite{CPRS2,CPRS3} we 
showed that if the spectral triple has a property we termed
 `isolated spectral dimension',
which is weaker than the assumption of `discrete dimension
spectrum' used in \cite{CM}, then we can evaluate 
the resolvent cocycle term-by-term by taking residues
at the critical point $r=(1-q)/2$. 
The resulting formula gives an index of a Breuer-Fredholm operator
 expressed in terms of the residues of zeta functions
at the critical point and
these residue functionals
assemble to give a version of the 
Connes-Moscovici residue cocycle.
We note that our formula is not identical with that
of \cite{CM} because we need to deal, in the von Neumann context,
with the problem that zero may be in the continuous spectrum of $\D$.

In this paper there are two main
theorems.  The first  proves that the resolvent cocycle
`almost' represents $(r-(1-q)/2)^{-1}$ times the Chern character $Ch$
of our semifinite spectral triples
in the sense  that they are cohomologous modulo functions 
which are holomorphic at 
$r=(1-q)/2$. We use the notation
$$(\phi^r_m)_{m=0}^M\sim \left(\frac{1}{r-(1-q)/2}Ch\right)$$ 
to represent this fact.

One consequence of our results is to provide an alternative 
proof of the semifinite local index formula in the case of odd spectral 
triples. This is because in \cite{CPRS2}, in order to show that the 
residue cocycle we obtain from the resolvent cocycle calculates the 
appropriate Breuer-Fredholm index
we had to start with the spectral flow formula of \cite{CP2}. This
spectral flow formula is quite 
difficult to prove in the case of general semifinite spectral
triples and it is desirable to find a more direct argument. 
Our second major theorem does this by proving 
that the residue cocycle represents the Chern character
so we may use the known fact
\cite{Co1} that the Chern character is an index cocycle to bypass the 
formula of
\cite{CP2}. (We note that, although the argument in \cite{Co1} 
that shows that the Chern character calculates a Fredholm index,
is formulated only in the standard case of spectral triples with 
$\mathcal N$ being all bounded operators on a separable Hilbert space,
 the results
of \cite{CPRS3} enable one to see that the arguments of \cite{Co1} go 
through without essential change for semifinite spectral triples.)

There are several benefits from the detailed treatment we give 
the proofs of our main results.

$\bullet$ We fill in many technical details 
for the transgression to the Chern character
that are absent from \cite{hig}. 

$\bullet$ In a
 semifinite spectral triple we have a semifinite von Neumann algebra $\cn$ 
with $\ca\subset\cn$  and $(1+\CD^2)^{-1/2}$ compact in $\cn$ satisfying some
summability hypothesis. To 
cope with the resulting `zero in the spectrum problem' for $\CD$,
the form of the residue cocycle and the proof that it represents the Chern
character are significantly modified from the standard type I situation.

$\bullet$ As noted above, we demonstrate the 
relationship between the resolvent 
cocycle and Higson's `improper cocycle'.
 In the process we obtain a renormalised 
version of this `improper cocycle' (which we call the reduced resolvent 
cocycle) clarifying  
the connection between our point of view and that of Higson.

$\bullet$ We obtain analytic continuation results for the resolvent cocycle 
evaluated on Hochschild cycles (of top degree) and on $(b,B)$ cycles. These 
results do not require assumptions on analytic continuation properties of 
individual zeta functions.

The exposition is organised as follows. We put preliminary material,
notation
and definitions needed for the main results in Section 2
including a
brief outline of the pseudodifferential calculus. Those 
parts of cyclic theory we need are introduced in Section 3.
Those expert in these matters can move straight to Section 4 where we
introduce the various cocycles which arose in our new proof of the local 
index theorem. Then we state our main theorem on the residue cocycle.
The proof is in Section 5 where we first prove
that in both the even and odd cases,  
the resolvent cocycle `almost' represents the Chern character
and we show that this implies  
that the residue cocycle is in the class of
the Chern character.
Section 6 derives some corollaries including the relationship to Higson's point of view.

\section{Definitions and Background}
\subsection{Semifinite Spectral Triples} 


We begin with some semifinite versions of standard definitions and results.
Let ${\mathcal K}_{\mathcal N }$ be the 
$\tau$-compact operators in ${\mathcal N}$
(that is the norm closed ideal generated by the projections
$E\in\mathcal N$ with $\tau(E)<\infty$). Here $\tau$ is a fixed faithful, 
normal, semifinite trace on the von Neumann algebra ${\mathcal N}$.

\begin{defn} A {\bf semifinite
spectral triple} $(\A,\HH,\D)$ is given by a Hilbert space $\HH$, a 
$*$-algebra $\A\subset \cn$ where $(\cn,\tau)$ is a semifinite 
von Neumann algebra with trace $\tau$ acting on
$\HH$, and a densely defined unbounded self-adjoint operator $\D$ affiliated 
to $\cn$ such that

1) $[\D,a]$ is densely defined and extends to a bounded operator in $\cn$ for 
all $a\in\A$

2) $(\lambda-\D)^{-1}\in\K_\cn$ for all $\lambda\not\in{\R}$

3) The triple is said to be {\bf even} if there is $\gamma\in\cn$ such that 
$\gamma^*=\gamma$, $\gamma^2=1$,  $a\gamma=\gamma a$ for all $a\in\A$ and 
$\D\gamma+\gamma\D=0$. Otherwise it is {\bf odd}.
\end{defn}

\subsection{Notes and  Remarks}. 
Henceforth we omit the term semifinite
as it is implied by the use of a faithful normal semifinite trace $\tau$
on $\cn$ in all of the subsequent text. In this paper, for simplicity of 
exposition,
we will deal only with unital algebras
${\mathcal A}\subset \cn$ where the identity of $\mathcal A$ is that of $\cn$.


\begin{defn}\label{qck} A semifinite spectral triple $(\A,\HH,\D)$ is 
$\bf{QC^k}$ for 
$k\geq 1$ 
($Q$ for quantum) if for all $a\in\A$ 
the operators $a$ and $[\D,a]$ are in the domain of $\delta^k$, where 
$\delta(T)=[\dd,T]$ is the partial derivation on $\cn$ defined by $\dd$. We say 
that 
$(\A,\HH,\D)$ is $\bf{QC^\infty}$ if it is $QC^k$ for all $k\geq 1$.
\end{defn}

{\bf Note}. The notation is meant to be analogous to the classical case, but 
we introduce 
the $Q$ so that there is no confusion between quantum differentiability of 
$a\in\A$ and classical differentiability of functions. We will sometimes
use the notation $da:=[\mathcal D,a]$ for the derivation $[\mathcal D,\cdot].$

\noindent{\bf Remarks concerning derivations and commutators}.  By partial 
derivation we mean that $\delta$ is defined on some subalgebra of $\cn$ which 
need not be (weakly) dense in $\cn$. More precisely, 
$\mbox{dom}\delta=\{T\in\cn:\delta(T)\mbox{ is bounded}\}$. We 
also note that if $T\in{\mathcal N}$, one can show that $[\dd,T]$ is bounded 
if and only if $[(1+\D^2)^{1/2},T]$ is bounded, by using the functional 
calculus to show that $\dd-(1+\D^2)^{1/2}$ extends to a bounded operator in 
$\cn$. In fact, writing $\dd_1=(1+\D^2)^{1/2}$ and $\delta_1(T)=[\dd_1,T]$ we 
have $\mbox{dom}\delta^n=\mbox{dom}\delta_1^n\ \ \ \ \forall n.$


Thus the condition defining $QC^\infty$ can be replaced by
\ben a,[\D,a]\in\bigcap_{n\geq 0}\mbox{dom}\delta_1^n\ \ \ \forall a\in\A.\een
This is important in situations where we cannot assume $\dd$ is invertible.

We also observe that if $T\in\cn$ and $[\D,T]$ is bounded, then $[\D,T]\in\cn$. 
Similar comments apply to $[\dd,T]$, $[(1+\D^2)^{1/2},T]$ and the more exotic 
combinations such as $[\D^2,T](1+\D^2)^{-1/2}$ which we will encounter later. 
The proofs of these statements can be found in \cite{CPRS2}.

Recall from \cite{FK} that if $S\in\mathcal N$, the {\bf t-th generalized
singular value} of S for each real $t>0$ is given by
$$\mu_t(S)=\inf\{||SE||\ \vert \ E \mbox{ is a projection in }
{\mathcal N} \mbox { with } \tau(1-E)\leq t\}.$$

The ideal $\LL^1({\mathcal N})$ consists of those 
operators $T\in {\mathcal N}$ such that $\n T\n_1:=\tau( |T|)<\infty$
where $|T|=\sqrt{T^*T}$.
In the Type I setting this is the usual trace class ideal. We will 
simply write $\LL^1$ for this ideal in order to simplify the notation, and 
denote the norm on $\LL^1$ by $\n\cdot\n_1$. An alternative definition in terms 
of singular values is that $T\in\LL^1$ if 
$\|T\|_1:=\int_0^\infty \mu_t(T) dt <\infty.$

Note that in the case where 
${\mathcal N}\neq{\mathcal B}({\mathcal H})$, 
$\LL^1$ need not be 
complete in this norm but it is complete in the norm $\|.\|_1 + \|.\|_\infty$.
(where $\|.\|_\infty$ is the uniform norm).

\subsection{The Pseudodifferential Calculus}

We refer to \cite[Section 6]{CPRS2} and \cite{CM,C4} for a full discussion of 
the pseudodifferential calculus, but present a brief review of the central 
ideas here. 
Given a densely-defined self-adjoint unbounded operator $D$ on a Hilbert space 
$\HH$, we set $\HH_\infty=\cap_{n\geq 0}\mbox{dom}D^n$. We denote by 
$\delta_1$ the derivation given by $T\mapsto [(1+D^2)^{1/2},T]$. As in
\cite[Section 6]{CPRS2} we let $|D|_1=(1+D^2)^{1/2},$ so that
$\delta_1(T)=[|D|_1,T].$ We then define 
linear spaces of operators for $r\in{\R}$
\ben OP^0=\cap_{n\geq 0}\mbox{dom}\delta_1^n,\ \ \ \ OP^r=(1+D^2)^{r/2}OP^0.
\een
We observe that $\delta_1$ clearly leaves $OP^0$ invariant.
If $(\A,\HH,\D)$ is a $QC^\infty$ spectral triple and we define $\HH_\infty$ 
and the spaces $OP^r$ using $\D$, then $\A,[\D,\A]\subset OP^0$. 
Defining $\nabla(T)=[\D^2,T],$ and setting $T^{(n)}=\nabla^n(T),$
we find that $a^{(n)}$ and $[\D,a]^{(n)}$ are in $OP^n$ for all $a\in\A$.

We recall Lemma 6.2 of \cite{CPRS2}.

\begin{lemma}[compare Lemma 1.1 of \cite{C4}] Let 
$b\in OP^0$. With $\s_1(b)=|D|_1b|D|_1^{-1}$ 
and $\varepsilon_1(b)=\delta_1(b)|D|_1^{-1}$ we have

1) $\s_1=Id+\varepsilon_1$,

2) $\varepsilon_1^n(b)=\delta_1^n(b)|D|_1^{-n}\in OP^0\ \ \ \forall n$,

3) $\s_1^n(b)=(Id+\varepsilon_1)^n(b)=
\sum_{k=0}^n{n\choose k}\delta_1^k(b)
|D|_1^{-k}\in OP^0\ \ \ \forall n.$
\end{lemma}

In a similar manner we can also prove the following.

\begin{lemma}
Let $b\in OP^0$. With $\gamma_1(b)=2\delta_1(b)|D|_1$, we have:

1) $\nabla=\delta_1^2 + \gamma_1$,

2) $\gamma_1^k(b)=2^k\delta_1^k(b)|D|_1^k \in OP^k\ \ \ \forall k$,

3) $b^{(n)}:=\nabla^n(b)=\left(\sum_{k=0}^n 2^k
{n\choose k}\delta_1^{2n-k}(b)
|D|_1^{-n+k}\right)|D|_1^n \in OP^n\ \ \ \forall n.$
\end{lemma}

\begin{proof}
Item 1) is a straightforward calculation noting that $\nabla(b)=[|D|_1^2,b].$
Item 2) follows from the definition of $\gamma_1$ by induction. Item 3) follows
from applying the binomial theorem to 1) and then using 2).
\end{proof}

\begin{defn}\label{deltatop}
We define an increasing sequence of norms $\|\cdot\|_k$ on $OP^0$ via
$\|b\|_k = \sum_{j=0}^k\|\delta_1^j(b)\|$ for $k\geq 0.$ This is closely
related to the  $\boldsymbol{\delta}${\bf-topology} on $QC^\infty$ algebras, 
$\A$, given 
by the family of seminorms:
$$a\to\n\delta^k(a)\n\ \ \ a\to\n\delta^k([\D,a])\n,\ \ \ k=0,1,2,\dots$$
\end{defn}

\begin{rem*} For $b\in OP^0$, we let 
$b_{(n)}=\sum_{k=0}^n 2^k{n\choose k}\delta_1^{2n-k}(b)|D|_1^{-n+k}$ so that
$b^{(n)}=b_{(n)}|D|_1^n$ where $b_{(n)}\in OP^0,$ and $||b_{(n)}||\leq C_n
||b||_{2n}$ where the constant $C_n$ depends only on $n.$
\end{rem*}

\begin{lemma}
Let $a\in OP^0$, and $n,p \geq 0$. Then

1) $\sigma_1^p(a_{(n)}) =\sum_{k=0}^n\sum_{j=0}^p 2^k{n\choose k}{p\choose j}
\delta_1^{2n-k+j}(a)|D|_1^{-n+k-j}$ is in $OP^0$, and

2) there is a positive constant $C_{p,n}$ depending only on $p$ and $n$
so that $\|\sigma_1^p(a_{(n)})\|\leq C_{p,n} \|a\|_{2n+p}.$
\end{lemma}

\begin{proof}
Item 1) is just a calculation combining the previous two lemmas and noting 
that $\sigma_1$ and $\delta_1$ not only commute but are both 
$|D|_1^{-1}$-linear. Item 2) follows from item 1) and the fact that 
$\| |D|_1^{-1}\|\leq 1.$
\end{proof}

\begin{cor}
If $b_0,b_1,...,b_m$ are in $OP^0$ and $n_0,n_1,...,n_m$ are nonegative
integers with \\$|n|:=n_0+n_1+\cdots+n_m,$ then there is a $C>0$ depending
only on $m$ and $|n|$ so that

$b_0^{(n_0)}b_1^{(n_1)}\cdots b_m^{(n_m)}=B |D|_1^{|n|}$ where $B\in OP^0$
and $\|B\|\leq C\|b_0\|_{2|n|} \|b_1\|_{2|n|}\cdots\|b_m\|_{2|n|}.$ 
\end{cor}

\begin{proof}
In the notation of the previous remark
\begin{eqnarray*}
b_0^{(n_0)}b_1^{(n_1)}\cdots b_m^{(n_m)}&=&
(b_0)_{(n_0)}|D|_1^{n_0}(b_1)_{(n_1)}|D|_1^{n_1}\cdots (b_m)_{(n_m)}|D|_1^{n_m}\\
&=&(b_0)_{(n_0)}\sigma_1^{n_0}((b_1)_{(n_1)})\sigma_1^{n_0+n_1}((b_2)_{(n_2)})
\cdots\sigma_1^{|n|}((b_m)_{(n_m)})|D|_1^{|n|}.
\end{eqnarray*}
The result now follows from the previous lemma
with 
$$B=(b_0)_{(n_0)}\sigma_1^{n_0}((b_1)_{(n_1)})
\sigma_1^{n_0+n_1}((b_2)_{(n_2)})\cdots\sigma_1^{|n|}((b_m)_{(n_m)}),$$
since $\|\cdot\|_{2k}\leq \|\cdot\|_{2|n|}$
for each $k\leq |n|.$
\end{proof}

A slight variation on the previous corollary is the following.

\begin{cor}
If $b_0,b_1,...,b_m$ are in $OP^0$ and $n_0,n_1,...,n_m$ are nonegative
integers with \\$|n|:=n_0+n_1+\cdots+n_m,$ then there is a $C>0$ depending
only on $m$ and $|n|$ so that
$$b_0^{(n_0)}b_1^{(n_1)}\cdots b_k^{(n_k)}D b_{k+1}^{(n_{k+1})}\cdots
b_m^{(n_m)}=B |D|^{|n|+1}$$
where $B\in OP^0$ and 
$$\|B\|\leq C\|b_0\|_{2|n|+1} \|b_1\|_{2|n|+1}\cdots\|b_m\|_{2|n|+1}.$$ 
\end{cor}
\begin{proof}
Write $D=F_0|D|_1$ where $F_0:=D|D|_1^{-1}\in OP^0$ and proceed as in the 
previous proof.
\end{proof}
If $\D$ is $n$-summable, so $(\lambda-\D)^{-1}\in\LL^n(\cn)$ for all 
$\lambda\not\in{\R}$, then any $T\in OP^r$ is $n/r$-summable.

\section{Cyclic Cohomology and Chern Characters}\label{cycliccohomology}

A major feature of \cite{Co4} 
is the  association to a suitable 
representative 
of a $K$-theory class, respectively a $K$-homology class, a class in 
{\bf periodic} cyclic homology, 
respectively a class in {\bf periodic} cyclic cohomology, 
called a {\bf Chern character} in both cases. The principal result is then  
\be \langle [x],[(\A,\HH,\D)]\rangle=
\langle [Ch_*(x)],[Ch^*(\A,\HH,\D)]\rangle,\label{indpair}\ee
where $[x]\in K_*(\A)$ is a $K$-theory class with representative $x$ and 
$[(\A,\HH,\D)]$ is the $K$-homology class of the spectral triple $(\A,\HH,\D)$. 
(The exact normalisations for these pairings depends on what kind of cochains 
one uses to represent cyclic cohomology.)
 
On the right hand side, $Ch_*(x)$ is the Chern character of $x$, and $[Ch_*(x)]$
its periodic cyclic 
homology class. Similarly $[Ch^*(\A,\HH,\D)]$ is the periodic cyclic cohomology 
class of 
the 
Chern 
character of $(\A,\HH,\D)$. 

We will describe the complexes defining cyclic cohomology that we use below. 
These are the cyclic complex and the $(b,B)$ bicomplex. We  also describe  the 
Chern character.

To define the (normalised) $(b,B)$ bicomplex, we introduce the following
linear spaces, \cite{Lo}.
Let $C_m=\A\otimes \bar\A^{\otimes m}$
where $\bar\A$ is the quotient $\A/\C\!\cdot\!\! I$ with $I$ being the identity
element of $\A$ and, assuming with no loss of generality that $\A$ is complete
in the  $\delta$-topology, \cite{R}, we employ the projective tensor product.
Let $C^m=Hom(C_m,\C)$ be the linear space of continuous linear functionals
><DEFANGED.834 on $C_m$.
We may define the $(b,B)$ bicomplex using these spaces (as opposed to $C_m=\A^{\otimes m+1}$ et cetera) and the resulting cohomology
will be the same. This follows because the bicomplex defined using
$\A\otimes \bar\A^{\otimes m}$ is quasi-isomorphic to that defined
using  $\A\otimes \A^{\otimes m}$ \cite{Lo}.
Similar comments apply to the cyclic complex.

We first define cyclic cohomology using the {\bf cyclic} complex.
A {\bf normalised cyclic cochain} on $\A$ is a functional $\psi\in C^m$ 
such that
\ben \psi(a_0,...,a_m)=(-1)^{m}\psi(a_m,a_0,...,a_{m-1}).\een
The set of all normalised cyclic cochains in $C^m$ is denote $C^m_{\lambda}$.
We say that $\psi$ is a {\bf cyclic cocycle} if for all $a_0,...,a_{m+1}\in\A$ 
we have 
$(b\psi)(a_0,...,a_{m+1})=0$ where
$$(b\psi)(a_0,a_1,\ldots,a_{m+1})=\hfill$$
$$\sum_{j=0}^{m}(-1)^j\psi(a_0,a_1,\ldots,a_ja_{j+1},\ldots,a_{m+1})
+(-1)^{m+1}\phi(a_{m+1}a_0,a_1,\ldots,a_{m}).$$
Then $H_\lambda^m(\A)$, the $m$-th cohomology group of $(C^m_\lambda,b)$
is defined to be the {\bf $m$-th cyclic cohomology group} of $\A$.

The cup product with the generator $\s\in H^2_\lambda({\C})$ of the cyclic 
cohomology of ${\C}$ defines a map $S:H^m_\lambda(\A)\to H^{m+2}_\lambda(\A)$ 
for any (locally convex) algebra $\A$. This {\bf periodicity operator} allows us 
to define the {\bf periodic cyclic cohomology} of $\A$ as the direct limit of 
the cyclic cohomology groups:
\ben H^*_{per}(\A)=\lim_{\to}(H^m_\lambda(\A),S),\een
where $*$ on the left hand side takes only the values even or odd, and on the 
right hand side we consider only those $m$ with the same parity as $*$.

A {\bf normalised $(b,B)$-cochain}, $\phi$ is a finite collection of 
multilinear functionals,
$$\phi=\{\phi_m\}_{m=0,1,...,M}\mbox{ with }\phi_m\in C^m.$$

An {\bf odd cochain} has $\phi_m=0$ for even $m$, while an 
{\bf even cochain} has $\phi_m=0$ for odd $m$.
It is a {\bf (normalised) $(b,B)$-cocycle} if, for all $m$,
$b\phi_m+B\phi_{m+2}=0$ where 
$b: C^m\to C^{m+1}$, $B:C^m\to C^{m-1}$
are the {\bf coboundary operators}. The operator $b$ is described above and $B$ 
is given by
$$(B\phi_m)(a_0,a_1,\ldots,a_{m-1})
=\sum_{j=0}^{m-1} 
(-1)^{(m-1)j}\phi_m(1,a_j,a_{j+1},\ldots,a_{m-1},a_0,\ldots,a_{j-1})$$
 We write $(b+B)\phi=0$ for brevity.
Thought of as functionals on $\A^{\otimes m+1}$
a normalised cocycle will satisfy
$\phi(a_0,a_1,\ldots,a_n)=0$ whenever any $a_j=1$ for $j\geq 1$.

 Similarly, a {\bf $(b^T,B^T)$-chain} $c$ is a (possibly infinite) collection 
$c=\{c_m\}_{m=0,1,...}$ with 
$c_m\in C_m$. 
The $(b,B)$-chain $\{c_m\}$ 
is a {\bf $(b^T,B^T)$-cycle} if $b^Tc_{m+2}+B^Tc_m=0$ for all $m$. More 
briefly, we write $(b^T+B^T)c=0$.
Here $b^T,B^T$ are the {\bf boundary operators} of cyclic homology, 
and are the transpose of the coboundary operators 
$b,B$ in the following sense. 

The {\bf pairing} between a $(b,B)$-cochain $\phi=\{\phi_m\}^M_{m=0}$ and a 
$(b^T,B^T)$-chain 
$c=\{c_m\}$ 
is given by
\ben \langle \phi,c\rangle= \sum_{m=0}^M\phi_m(c_m).\een
This pairing satisfies
\ben \langle (b+B)\phi,c\rangle=\langle\phi,(b^T+B^T)c\rangle.\een
All of the cocycles we consider in this paper are in fact defined
as functionals on $\oplus_m\A\otimes \bar\A^{\otimes m}$.
Henceforth we will drop the superscript on  $b^T,B^T$ and just write $b,B$
for both boundary and coboundary operators
as the meaning will be clear from the context.

Our next aim is to define the Chern character of a finitely summable 
Fredholm module.
First we need a definition.

\begin{defn}
A {\bf pre-Fredholm module} for a unital $*$-algebra
$\mathcal A$  is a pair $(\HH, F)$ where $\mathcal A$ is
represented in $\mathcal N$ (a semifinite von Neumann
algebra acting on $\HH$ with fixed trace $\tau$) and $F$ is a self-adjoint
Breuer-Fredholm operator in $\mathcal N$ satisfying:

$1.\: 1-F^2 \in {\mathcal {K_N}},\:and$

$2.\: [F,a] \in {\mathcal {K_N}}\: for\: a \in {\mathcal A}.$

\noindent We say that $(\HH, F)$ is {\bf even} if there is a grading operator
><DEFANGED.835 $\gamma\in\mathcal N$ such that $\gamma^*=\gamma$, $\gamma^2=1$,
$[\gamma,a]=0$ for all $a\in\mathcal A$ and $\gamma F+F\gamma=0$; otherwise, 
$(\HH, F)$ is {\bf odd}. Our formulas will often include a factor of
$\gamma$ in the odd case as well as the even case: in the odd case
we interpret $\gamma$ to be $1$. If $1-F^2=0$ we drop the prefix "pre-". 
If $[F,a]\in\LL^p(\NN)$ 
for all $a\in\A$, we say that $(\HH,F)$ is $\boldsymbol p${\bf-summable}.
\end{defn}

{\bf Pertinent Example}. Semifinite spectral triples give rise to pre-Fredholm 
modules via
\ben (\A,\HH,\D)\mapsto (\HH,F=\D(1+\D^2)^{-1/2}).\een
One views spectral triples as geometric representatives of $K$-homology 
classes, 
in much the same way that one views differential forms as geometric 
representatives of cohomology classes.

If the semifinite spectral triple $(\A,\HH,\D)$ is $QC^\infty$ and finitely 
summable with $(1+\D^2)^{-s/2}$ trace class for all $s>q$, and has $\D$ 
invertible, then
\ben (\HH, F=\D\dd^{-1})\een
is a $[q]+1$-summable Fredholm module where $[\cdot]$ denotes the integer part. 

\begin{defn} \label{conditional} We define the `conditional trace'  $\tau'$ by
\ben \tau'(T)=\frac{1}{2}\tau(F(FT+TF)),\een
provided $FT+TF\in\LL^1(\NN)$ (as it will be in our case, see 
\cite[p293]{Co4}). 
Note that if $T\in\LL^1(\NN)$ 
we have (using the trace property and $F^2=1$)
\be \tau'(T)=\tau(T).\label{tracesequal}\ee
\end{defn}

The {\bf Chern character} $[Ch_F]$ of an $n+1$-summable Fredholm module 
$(\HH,F)$ 
($n$ an integer) is 
the class in periodic cyclic cohomology of 
the cyclic cocycles
\ben \lambda_m\tau'(\gamma a_0[F,a_1]\cdots[F,a_m]), \ \ \ a_0,...,a_m
\in\A,\ \  m\geq n,\ \  m\mbox{ even if }(\HH,F)\mbox{ even, and odd 
otherwise}.\een

Here $\lambda_m$ are constants ensuring that this collection of cocycles 
yields a well-defined periodic class, and they are given by
\ben \lambda_m=\left\{\begin{array}{ll} (-1)^{m(m-1)/2}\Gamma(\frac{m}{2}+1) 
& \ \ m\ \ \ {\rm even}\\ \sqrt{2i}(-1)^{m(m-1)/2}\Gamma(\frac{m}{2}+1) &
\ \ m\ \ \ {\rm odd}\end{array}\right..\een
The {\bf class of the Chern character} of an $n+1$-summable Fredholm module is 
represented by 
the cyclic cocycle in bottom dimension $n,$ $Ch_F\in C^{n}_\lambda(\A)$ :
\ben Ch_F(a_0,...,a_{n})=\lambda_{n}\tau'(\gamma a_0[F,a_1]
\cdots[F,a_n]), \ \ \ \ \ a_0,...,a_n\in\A.\een
We will always take the cyclic cochain $Ch_F$ (or its $(b,B)$ analogue; see 
below) as representative of $[Ch_F]$, and will often refer to $Ch_F$ as the 
Chern character.

Since the Chern character is a cyclic cochain, it lies in the image of the 
operator $B$, \cite[Corollary 20, III.1.$\beta$]{Co4}, and so $BCh_F=0$ since 
$B^2=0$. Since $bCh_F=0$, we may regard the Chern character as a one term 
element of the $(b,B)$ bicomplex. However, 
the correct normalisation is (taking the Chern character to be in degree $n$)
\ben C^n_\lambda\ni Ch_F\mapsto \frac{(-1)^{[n/2]}}{n!}Ch_F\in C^n.\een
Thus instead of $\lambda_n$ defined above, we use $\mu_n$
\ben \mu_n=\frac{(-1)^{[n/2]}}{n!}\lambda_n=\left\{\begin{array}{ll} 
\frac{\Gamma(\frac{n}{2}+1)}{n!} 
& \ \ n\ \ \ {\rm even}\\ & \\ \sqrt{2i}\frac{\Gamma(\frac{n}{2}+1)}{n!} &
\ \ n\ \ \ {\rm odd}\end{array}\right..\een
The difference in normalisation between periodic and $(b,B)$ is due  to 
the way the index pairing is defined in the two cases, \cite{Co4}, and
compatibility with the periodicity operator.

Our next task is to show that if our spectral triple $(\A,\HH,\D)$
is such that $\D$ is not invertible, we can replace it by a new spectral 
triple in the 
same $K$-homology class in which the unbounded operator is invertible.
This is not a precise statement in the general semifinite case, as our 
spectral triples will not define $K$-homology classes in the usual sense. When 
we say that two spectral triples are in the same $K$-homology class, we shall 
take this to mean that the associated pre-Fredholm modules are operator 
homotopic up to the addition of degenerate Fredholm modules (see \cite{KK} 
for these notions, which make sense in our context).

\begin{defn} Let $(\A,\HH,\D)$ be a spectral triple. For any 
$\mu\in{\R}\setminus\{0\}$,  define the `double' of 
$(\A,\HH,\D)$ to be the spectral triple $(\A,\HH^2,\D_\mu)$ with 
$\HH^2=\HH\oplus\HH$, 
and the action of $\A$ and $\D_\mu$ given by
\ben \D_\mu=\bma \D & \mu\\ \mu & -\D\ema,\ \ \ \ a\mapsto\bma a & 0\\ 0 & 0\ema,
\ \ \forall a\in\A.\een
\end{defn}
{\bf Remark} Whether $\D$ is invertible or not, $\D_\mu$ always is invertible, 
and $F_\mu=\D_\mu|\D_\mu|^{-1}$
has square 1. This is the chief reason for introducing this construction. 
\begin{lemma}\label{noninv} The $K$-homology classes of $(\A,\HH,\D)$ and 
$(\A,\HH^2,\D_\mu)$ are the same. A representative of this class is 
$(\HH^2,F_\mu)$ 
with $F_\mu=\D_\mu|\D_\mu|^{-1}$.
\end{lemma}



The most basic consequence of Lemma \ref{noninv}, whose proof may be found in 
\cite{CPRS1}, comes from the following (see \cite[pp65-68]{Co1}
for the proof).

\begin{prop}\label{chind} The periodic cyclic cohomology class of the 
Chern character of a finitely summable Fredholm module depends only on its 
$K$-homology class.
\end{prop}

In particular, therefore,
the Chern characters of $(\A,\HH,\D)$ and $(\A,\HH^2,\D_\mu)$ have the same 
class in periodic cyclic cohomology, and this can be computed (indeed is 
defined!) using the 
Fredholm module $(\HH^2,F_\mu)$, and this class is independent of $\mu$. 

\section{The theorem on the residue cocycle}

\subsection{The odd semifinite local index formula}
Our main results in this paper
are motivated by, and have consequences for,
 the semifinite local index formula in the odd case \cite{CPRS2}.

To state it we need some notation, beginning with multi-indices $(k_1,...,k_m)$, 
$k_i=0,1,2,...$, whose length $m$ will always be
clear from the context. Next we
 write $|k|=k_1+\cdots+k_m$, and define $\alpha(k)$ by
\ben \alpha(k)=\frac{1}{k_1!k_2!\cdots k_m!(k_1+1)(k_1+k_2+2)\cdots(|k|+m)}.\een
The numbers $\s_{n,j}$ are defined by the equality
\ben 
\prod_{j=0}^{n-1}(z+j+1/2)=\sum_{j=0}^{n}z^j\s_{n,j}.\een

If $(\A,\HH,\D)$ is a $QC^\infty$ spectral triple and $T\in\cn$, we write 
$T^{(n)}$ to denote  the iterated commutator 
$[\D^2,[\D^2,[\cdots,[\D^2,T]\cdots]]]$ where we have $n$ commutators with 
$\D^2$. It follows from the remarks after Definition \ref{qck} that operators of 
the form $T_1^{(n_1)}\cdots T_k^{(n_k)}(1+\D^2)^{-(n_1+\cdots+n_k)/2}$ are in 
$\cn$ for $T_i=[\D,a_i]$, $a_i\in\A$. 

\begin{defn}\label{dimension} 
If $(\A,\HH,\D)$ is a $QC^\infty$ spectral triple, 
we call 
\ben q=\inf\{k\in{\R}:\tau((1+\D^2)^{-k/2})<\infty\}\een
the {\bf spectral dimension} of $(\A,\HH,\D)$.
We say that $(\A,\HH,\D)$ has {\bf isolated spectral dimension} if
for $b$ of the form  
$$b=a_0[\D,a_1]^{(k_1)}\cdots[\D,a_m]^{(k_m)}(1+\D^2)^{-m/2-|k|}$$
the zeta functions
\ben \zeta_b(z-(1-q)/2)=\tau(b(1+\D^2)^{-z+(1-q)/2})\een
have analytic continuations to a deleted neighbourhood of $z=(1-q)/2$.
\end{defn}

{\bf Remark}. Observe that we allow the possibility that the analytic 
continuations of these zeta functions may have an essential singularity at 
$z=(1-q)/2$. All that is necessary for us is that the residues at this point 
exist. Note that discrete dimension spectrum implies isolated spectral 
dimension.

\begin{defn}\label{resnotation}
Now we define, for $(\A,\HH,\D)$ having isolated spectral dimension
and 
$$b=a_0[\D,a_1]^{(k_1)}\cdots[\D,a_m]^{(k_m)}(1+\D^2)^{-m/2-|k|},$$
the numbers, $\boldsymbol{\tau_j(b)}$ by:
\ben \tau_j(b)=res_{z=(1-q)/2}(z-(1-q)/2)^j\zeta_b(z-(1-q)/2).\een
\end{defn}
The hypothesis of isolated spectral dimension is clearly necessary here 
in order to define the residues.
Let $Q$ be the spectral projection of $\D$ corresponding to the
interval $[0,\infty)$.

In \cite{CPRS2} we proved the following result:

\begin{thm}[Odd semifinite local index formula]\label{SFLIT} 
Let $(\A,\HH,\D)$ be an 
odd finitely summable $QC^\infty$ spectral triple with spectral
dimension $q\geq 1$. Let $N=[q/2]+1$ where $[\cdot]$ denotes the 
integer part, and let $u\in\A$ be unitary. Then

1) \qquad\qquad\qquad $index(QuQ)=\frac{1}{\sqrt{2\pi i}} 
res_{r=(1-q)/2}\left( 
\sum_{m=1,odd}^{2N-1} \phi_m^r(Ch_m(u))\right)$

where for $a_0,...,a_m\in\A$, $l=\{a+iv:v\in{\R}\}$, $0<a<1/2$, 
$R_s(\lambda)=(\lambda-(1+s^2+\D^2))^{-1}$ and $r>0$ we define 
$\phi_m^r(a_0,a_1,...,a_m)$ to be
\ben \frac{-2\sqrt{2\pi i}}{\Gamma((m+1)/2)}
\int_0^\infty s^m\tau\left(\frac{1}{2\pi i}
\int_l\lambda^{-q/2-r}a_0R_s(\lambda)[\D,a_1]R_s(\lambda)\cdots
[\D,a_m]R_s(\lambda)d\lambda\right)ds.\een
In particular the sum on the right hand side of $1)$ analytically continues 
to a deleted neighbourhood of $r=(1-q)/2$ with {\em at worst} a simple pole 
at $r=(1-q)/2$.
Moreover, the complex function-valued cochain $(\phi_m^r)_{m=1,odd}^{2N-1}$ 
is a $(b,B)$ 
cocycle for $\A$ modulo functions holomorphic in a half-plane containing 
$r=(1-q)/2$.

2) The index is also the residue of a sum of zeta 
functions:
\bean &&\frac{1}{\sqrt{2\pi i}} res_{r=(1-q)/2} \left(\sum_{m=1,odd}^{2N-1}
\sum_{|k|=0}^{2N-1-m}\sum_{j=0}^{|k|+(m-1)/2}(-1)^{|k|+m}
\alpha(k)\Gamma((m+1)/2)\s_{|k|+(m-1)/2,j} \right.\nno
&& \qquad\qquad\biggl.(r-(1-q)/2)^j\tau\left(u^*[\D,u]^{(k_1)}
[\D,u^*]^{(k_2)}\cdots[\D,u]^{(k_m)}(1+\D^2)^{-m/2-|k|-r+(1-q)/2}\right)
\Biggr).\eean
In particular the sum of zeta functions on the right hand side analytically 
continues to a deleted neighbourhood of $r=(1-q)/2$ and has {\rm at worst} a 
simple pole at 
$r=(1-q)/2$.

3) If $(\A,\HH,\D)$ also has isolated spectral dimension then
\ben index(QuQ)=\frac{1}{\sqrt{2\pi i}}\sum_m \phi_m(Ch_m(u))\een
where for $a_0,...,a_m\in\A$
\bean \phi_m(a_0,...,a_m)&=&res_{r=(1-q)/2}\phi^r_m(a_0,...,a_m)=\sqrt{2\pi i}
\sum_{|k|=0}^{2N-1-m}(-1)^{|k|}
\alpha(k)\times\nno
&\times&\sum_{j=0}^{|k|+(m-1)/2}\s_{(|k|+(m-1)/2),j}\tau_j
\left(a_0[\D,a_1]^{(k_
1)}\cdots[\D,a_m]^{(k_m)}(1+\D^2)^{-|k|-m/2}\right),\eean
and $(\phi_m)_{m=1,odd}^{2N-1}$ is a $(b,B)$ cocycle for $\A$. When $q<2N-1$, 
the term with $m=2N-1$ is zero, and 
for $m=1,3,...,2N-3$, all the top terms with $|k|=2N-1-m$ are zero. 
\end{thm}

In \cite{CPRS3} we established the even version of this theorem however we will
have no need to state the full result here.

\subsection{The resolvent cocycle in the general case}

The following definition establishes some notation 
needed to treat the even and 
odd cases on the same footing (this was not done in our earlier work).

Our philosophy on the resolvent cocycle is that it provides a substitute 
for the
JLO cocycle (the foundation of the original Connes-Moscovici local 
index theorem). The stronger version of the local index theorem that we
obtain arises from the more detailed information inherent in a cocycle that
is sensitive to the spectral dimension which the JLO cocycle is not. 
The form of our resolvent cocycle was deduced from the 
spectral flow formula of \cite{CP2} in the odd case and from a generalised
McKean-Singer formula in the even case \cite{CPRS3}. 

\begin{defn} Let $(\A,\HH,\D)$ be a spectral triple with spectral dimension 
$q\geq 1$. Let $P$ denote the 
{\bf parity} of the triple, so $P=0$ for even triples and $P=1$ for odd 
triples. 
Let $A$ denote $(P-1)$, the {\bf anti-parity} so $A=1$ for even triples 
and $A=0$ for odd triples.
We adopt the convention that $\dd$ and elements of $\A$ have {\bf grading 
degree} 
zero, while $\D$ has grading degree one. In the even case this is of course 
the actual grading degree of the spectral triple. We denote the grading degree 
of $T\in OP^*$ by $deg(T)$. Finally, let $N=[(q+1+P)/2]$ where $[\cdot]$ 
denotes the integer part. So, $M=2N-P$ is the greatest integer of parity $P$
in $q+1.$ In particular, if $q$ is an integer of parity $P$ then $M=2N-P=q.$
\end{defn}

The grading degree is used to define the graded commutator:
\ben [T,R]_\pm:=TR-(-1)^{deg(T)deg(R)}RT.\een
In particular, we have $[\D,a]_\pm =[\D,a]$ for all $a\in \A.
$ Similarly, since
$deg(\D^2)=0$, we have $[\D^2,T]_\pm = [\D^2,T]$ for all $T.$
The following definition generalises the expectations introduced in 
\cite{CPRS2,CPRS3} to deal with both the even and odd cases
in a uniform fashion.

\begin{defn}\label{expectation} Let $0<a<1/2$ and let $l$ be the vertical line 
$l=\{a+iv:v\in\R\}$.
For $m\geq 0$, $s\in[0,\infty)$ and operators $A_0,...,A_m$, $A_i\in OP^{k_i}$, 
with $k_0+\cdots+k_m-2m<2Re(r)$  define
\ben \la A_0,...,A_m\ra_{m,s,r}=\tau\left(\frac{1}{2\pi i}\gamma\int_l 
\lambda^{-q/2-r}
A_0R_s(\lambda)A_1\cdots A_mR_s(\lambda)d\lambda\right).\een
Here $\gamma$ is the ${\Z}_2$-grading in the even case and the identity operator 
in the odd case, and $R_s(\lambda)=(\lambda-(1+s^2+\D^2))^{-1}$. 
\end{defn}
We now state the definition of the resolvent cocycle in terms of the 
expectations $\la\cdots\ra_{m,s,r}$.

\begin{defn} (Compare Definition 7.9 \cite{CPRS2}) Let $(\A,\HH,\D)$ be a 
spectral triple with spectral dimension $q\geq 1$ and parity $P.$
Introduce constants $\eta_m, \ m=0,1,2,\ldots$ 
with  $P\equiv m(\bmod\ 2)$ by
$$ \eta_m=\left(-\sqrt{2i}\right)^P2^{m+1}\frac{\Gamma(m/2+1)}{\Gamma(m+1)}.$$
Then for $Re(r)>\frac{1}{2}(1-m),$ the $m$-th component of the 
{\bf resolvent cocycle} $\phi_m^r:\A^{\otimes m+1}\to{\C}$ is defined by: 
\ben \phi_m^r(a_0,...,a_m)=\eta_m\int_0^\infty s^m\la 
a_0,da_1,..,da_m\ra_{m,s,r}ds.\een
Where recall that $da:=[\D,a]$ for $a\in\A.$
\end{defn}

A basic result which is implicit in \cite{CPRS2,CPRS3} is that
$\phi^r_m$ is a continuous map from $\A^{\otimes m+1}$ with 
the $\delta$-topology (see Definition \ref{deltatop}) to the space of
functions $F_m$
defined and 
holomorphic in the half-plane $Re(r)>(1-m)/2$ (with the topology of
uniform convergence on compacta). This turns out to be surprisingly subtle
so we include the proof with a
slight generalisation which we will require a little later. Let 
$$R_{s,t}(\lambda):=(\lambda-(t+s^2+\D^2))^{-1}$$
where if $\D^2\geq\delta>0$ (and so invertible) we allow $t\in[0,1]$
and if $\D$ is not invertible we allow only $t=1$. Define
$\phi^r_{m,t}$ just as $\phi^r_m$, but using $R_{s,t}(\lambda)$ in
place of $R_s(\lambda)$, and similarly for $\la\dots\ra_{m,s,r,t}$.

\begin{lemma}\label{better-be-true} Let $OP^0$ have the $\delta$-topology.
(see Definition \ref{deltatop}) 
 Then with $t\in[0,1]$ as above, the map
$$ (OP^0)^{\otimes m+1}\ni (A_0,\dots,A_m)\mapsto \int_0^\infty s^m\la
A_0,\dots,A_m\ra_{m,s,r,t}ds \in F_m$$
is a continuous multilinear functional. In particular for each $r\in\C$ 
with $Re(r)>(1-m)/2$,\\ $\phi^r_{m,t}$
restricts to a continuous multilinear functional on $\A$ with the
$\delta$-topology. Moreover,\\ $(a_0,a_1,\dots,a_m)\mapsto(r\mapsto 
\phi^r_{m,t}(a_0,\dots, a_m))$ is an element of $Hom(\A^{\otimes m+1},F_m)$.
\end{lemma}
\begin{proof}

Since $\A,\ [\D,\A]\subset OP^0$, the penultimate statement of the Lemma 
follows easily from the first. 
To prove the first statement we begin by rewriting the expectations
$\la\dots\ra_{m,s,r,t}$ using both the $s$-trick (see Lemma \ref{strick} below),
and the $\lambda$-trick (see Lemma \ref{lambdatrick} below). 
The $s$-trick says that for all $A_0,\dots,A_m\in
OP^0$ we have
$$\int_0^\infty s^m\la A_0,\dots,A_m\ra_{m,s,r,t} ds
=\frac{-2}{m+1}\sum_{j=0}^m \int_0^\infty s^{m+2} \la
A_0,\dots,A_j,1,A_{j+1},\dots,A_m\ra_{m+1,s,r,t},$$
where both sides exist for $Re(r)>(1-m)/2$. On the other hand, the
$\lambda$-trick says that for all $A_0,\dots,A_m\in
OP^0$ we have
$$\int_0^\infty s^m\la A_0,\dots,A_m\ra_{m,s,r,t} ds
=\frac{-1}{q/2+r-1}\sum_{j=0}^m \int_0^\infty s^{m} \la
A_0,\dots,A_j,1,A_{j+1},\dots,A_m\ra_{m+1,s,r-1,t},$$
where both sides exist for $Re(r)>(1-m)/2$ provided the first factor on the right
does not introduce a pole in this region.

{\bf Some simple observations:}\\
(1) The $s$-trick and the $\lambda$-trick commute.\\
(2) Both tricks leave the region of convergence unaffected, with the
proviso that we don't introduce a pole with the $\lambda$-trick.\\
(3) After $X$ combined applications of the two tricks, the term with
any one particular pattern of $1's$ dispersed among the $A_i's$ will
appear with frequency exactly $X!$.
 
Let $M=2N-P$ so that $m$ also has parity $P$. If $M=m$ we do nothing to 
the expression for $\phi^r_{m,t}$. Assuming $m<M$ we will apply the $s$-trick 
exactly $(M-m)/2$ times (which will raise the power of $s^m$ to $s^M$ 
{\bf and no more}) and then
the $\lambda$-trick also exactly $(M-m)/2$ times,  which will raise our 
$m+1$-tuple to an $M+1$-tuple without introducing a pole in the region
$Re(r)>(1-m)/2$.

Thus, we get
by summing over all $k=(k_0,\dots,k_m)$ where each $0\leq k_i\leq (M-m)$ 
and $|k|=(M-m)$:
\begin{eqnarray*}
&&\int_0^\infty s^m \la A_0,\dots,A_m\ra_{m,s,r,t} ds\\
&&=(2)^{\frac{M-m}{2}}(M-m)!\prod_{b=1}^{\frac{M-m}{2}}\frac{1}{q/2+r-b}
\prod_{j=1}^{\frac{M-m}{2}}\frac{1}{m+j}\times\\
&&\times\sum_k\int_0^\infty s^M
\la A_0,1^{k_0},A_1,1^{k_1},\dots,A_m,1^{k_m}\ra_{M,s,r-(M-m)/2,t} ds
\end{eqnarray*}
Where, of course, we mean $1^{k_i}=1,1,\dots,1$ with $k_i$ one's. 
We observe that both sides converge for $Re(r)>(1-m)/2.$ Now the poles occur
when $r=b-q/2$ and since $b\leq(M-m)/2$ and $M\leq q+1$ we have for such poles,
$r\leq (M-m)/2 -q/2\leq (q+1-m)/2 -q/2 =(1-m)/2$ so that the pole is not
in the region, $Re(r)>(1-m)/2.$ That is, both sides exist in this region.

Now ignoring the prefactors, we have a sum of integrals
where we write $R$ for $R_{s,t}(\lambda)$ and each $n_i=k_i+1$
so that $n=(n_0,\dots,n_m)$ where $1\leq n_i\leq(M-m)+1$
and $|n|=M+1$:
$$\sum_n\int_0^\infty s^M\tau\left(\frac{1}{2\pi i}\gamma
\int_l\lambda^{-q/2-r+(M-m)/2}
A_0R^{n_0}A_1R^{n_1}\cdots A_mR^{n_m} d\lambda\right)ds$$

By Lemma 2.2 of \cite{CP1} we can estimate the trace of the $\lambda$-integral 
by integrating the trace norm of the integrand. 

For $n$ fixed let $p_{n_i}=(M+1)/n_i$ so that $\sum_i p_{n_i} =1$ 
and using H\"{o}lder's inequality we estimate:
\bean&&\sum_n\Vert A_0R^{n_1}A_1R^{n_2}\cdots
R^{n_m}A_mR^{n_{m+1}}\Vert_1\leq\sum_n\Vert A_0R^{n_1}\Vert_{p_{n_i}}\cdots
\Vert A_m R^{n_{m+1}}\Vert_{p_{n_{m+1}}}\nno
&\leq& \sum_n\Vert A_0\Vert\cdots\Vert A_m\Vert\cdot\Vert R^{n_1}\Vert_{p_{n_1}}
\cdots\Vert R^{n_{m+1}}\Vert_{p_{n_{m+1}}}= 
\sum_n\Vert A_0\Vert\cdots\Vert A_m\Vert\cdot\Vert R^{M+1}\Vert_1.\eean

Thus the iterated integral will be absolutely convergent if we show the 
boundedness of
\bean &&\int_0^\infty s^M\int_{-\infty}^\infty
\sqrt{a^2+v^2}^{-q/2-Re(r)+(M-m)/2}\Vert R^{M+1}\Vert_1dvds\nno
&\leq& 
\int_0^\infty s^M\int_{-\infty}^\infty
\sqrt{a^2+v^2}^{-q/2-Re(r)+(M-m)/2}
\sqrt{(s^2+a)^2+v^2}^{-M-1+q/2+\epsilon/2}dvds,\eean
the last inequality coming from Lemma \ref{lambda} below, and
$\epsilon>0$ is arbitrarily small. By Lemma \ref{intest} below, this integral is 
finite provided
$$M-2(M+1-q/2-\epsilon/2)<-1,\ \mbox{and}\
M-2(M+1-q/2-\epsilon/2)+1-2(q/2+Re(r)-(M-m)/2)<-2.$$
These conditions reduce to 
$$ q+\epsilon < M+1,\ \ \ (1-m)/2<Re(r)-\epsilon/2.$$
As $\epsilon>0$ is arbitrary, we see that the integral is finite for
$Re(r)>(1-m)/2$. 

With a little more effort using the above estimate and the ideas in the proof
of Lemma 5.4 in \cite{CPRS2} one can show that the double integral is 
{\bf uniformly bounded} in any closed right half-plane in the region
$Re(r)>(1-m)/2$ and is $O(1/Re(r))$ as $Re(r)\to\infty.$ That is,
$$|\phi^r_{m,t}(a_0,a_1,\dots,a_m)|\leq C(r_0)\n a_0\n\cdot\n da_1\n
\cdots\n da_m\n$$
uniformly in $r$ with $Re(r)\geq r_0>(1-m)/2$. Moreover, $C(r_0)$ is
$O(1/r_0)$ ar $r_0\to\infty$.

The fact that the map $r\mapsto \phi^r_{m,t}(a_0,\dots, a_m))$ is 
holomorphic in the
half-plane $Re(r)>(1-m)/2$ is similar to the proof of Lemma 7.4 of 
\cite{CPRS2}.
\end{proof}
\subsection{The residue cocycle}

The semifinite local index formula in the form  described  above
entails the introduction of the residue
cocycle, which we now describe.

\begin{defn} Let $(\A,\HH,\D)$ be a $QC^\infty$ finitely summable spectral 
triple with isolated spectral dimension $q\geq 1$. Let $M=2N-P$.
For $m=P,P+2,\dots,M$, and using the notation of Definition 
\ref{resnotation} define functionals 
$\phi_m$ by 
\ben \phi_m(a_0,...,a_m)=\sqrt{2\pi 
i}\sum_{|k|=0}^{M-m}\!\!(-1)^{|k|}\alpha(k)\!
\sum_{j=A}^{h}\s_{h,j}\tau_{j-A}\left(\gamma a_0[\D,a_1]^{(k_1)}\cdots[\D,a_m]
^{(k_m)}(1+\D^2)^{-|k|-m/2}\right),\een
where $h=|k|+(m-P)/2$. Here $\gamma$ denotes the $\Z_2$-grading in the even case 
and the identity operator in the odd case. Note that $M$ is the greatest odd
(respectively, even) integer in $(q+1)$ when the spectral triple is odd 
(respectively, even).


It follows from the results of \cite{CPRS2,CPRS3} that
 $\phi=(\phi_m)_{m=P}^M$ is a $(b,B)$-cocycle, called the {\bf residue cocycle}.
\end{defn}

The relationship between the resolvent cocycle and the residue cocycle 
is captured by the following result proved in \cite{CPRS2,CPRS3}.

\begin{thm} Let $(\A,\HH,\D)$ be a $QC^\infty$ finitely summable spectral triple 
with isolated spectral dimension $q\geq 1$. 
When evaluated on any $a_0,...,a_m\in\A$, with $m\equiv P(\bmod\ 2)$
 the components 
$\phi^r_m$ of the resolvent cocycle $(\phi^r)$ analytically continue to 
a deleted neighbourhood of
$r=(1-q)/2$.
 Moreover, if we denote this continuation by $\varphi^r_m(a_0,...,a_m)$ then
\ben res_{r=(1-q)/2}\varphi^r_m(a_0,...,a_m)
=\phi_m(a_0,...,a_m).\een
\end{thm}

It is the proof this last result that shows that the resolvent 
cocycle is indeed playing the same role as JLO does for the original proof of 
the local index theorem.  

\subsection{Statement of the theorem}

Our main result on the residue cocycle
establishes the following equality in cyclic cohomology.

\begin{thm}\label{res=chern} Let $(\A,\HH,\D)$ be a 
 $QC^\infty$ spectral triple (even or odd) with
spectral dimension $q\geq 1$ and
 isolated spectral dimension. 
Then the residue cocycle represents the Chern character of the $K$-homology 
class of $(\A,\HH,\D)$.
\end{thm}

\section{The residue and resolvent cocycles represent the Chern character}
\label{chernsection}

\subsection{Preamble}
Our methods are inspired by Higson, \cite{hig}, and 
we follow his approach quite closely. We present our arguments in full 
because numerous algebraic identities are, while very similar to Higson's, 
different in small details. Moreover $\D$ may have zero in its continuous
spectrum and this forces us to modify the standard approach. Finally
because we are building on
a proof of a version of  the Connes-Moscovici theorem valid under 
weaker hypotheses and of 
much greater generality we must work very hard to establish
some essential subtle estimates. 

The estimates of \cite{CPRS2} show that the
resolvent cocycle is a cocycle with values in the functions
defined and holomorphic on the right half plane $Re(r)>0$ (odd case) 
or $Re(r)>1/2$ (even case). In fact, for the individual 
functionals in the expression for the resolvent cocycle of degree $m,$
the half-plane of holomorphy increases as $m$ increases. We  will find
that our transgression arguments require refinements of the estimates
of \cite{CPRS2} whose proofs force us into a rather lengthy treatment.

Our immediate aim is to prove the following two statements
and from these we will deduce our theorem on the residue cocycle.
\begin{thm}\label{approx} Let $(\A,\HH,\D)$ be a $QC^\infty$ finitely summable 
spectral triple with dimension $q\geq 1$ and $\D$ {\bf invertible}. Let $M=2N-P$ 
where 
$N=[(q+1+P)/2],$ and note that in the odd (even) case, $M$ is the greatest
odd (even) integer in $q+1.$ Then in the $(b,B)$ bicomplex for $\A$ 
with coefficients in 
functions holomorphic for $Re(r)>1/2$, the resolvent cocycle 
$(\phi_m^r)_{m=P}^M$ is cohomologous to
\ben \frac{1}{(r-(1-q)/2)}Ch^M_F\een
modulo cochains with values in the functions holomorphic in a half-plane 
containing $r=(1-q)/2$. Here $F=\D\dd^{-1}$, and $Ch_F^M$ denotes the 
representative of the Chern character in dimension $M$.
\end{thm}

\begin{thm}\label{exactly} If $(\A,\HH,\D)$ is a $QC^\infty$ spectral triple 
with isolated spectral dimension $q\geq 1$ and $\D$ {\bf invertible}, then 
the cyclic 
cohomology class of the residue cocycle coincides with the class of the Chern 
character of $(\HH,F=\D\dd^{-1})$.
\end{thm}

There are two main steps involved in proving these statements. First 
in subsections 5.3 and 5.4 we need to 
define a `transgression' cochain which provides a cohomology (modulo cochains 
with values in the functions holomorphic in a half-plane containing 
$(1-q)/2$) between the 
resolvent cocycle and a single term cyclic cochain which is `almost' a cocycle. 
Second in subsection 5.5 we deform the resulting single term cyclic cochain 
to the Chern character. 
In this process we introduce error terms that are holomorphic at $r=(1-q)/2$.
Theorem \ref{exactly} follows
on taking residues and this requires the isolated 
spectral dimension hypothesis.

Both these steps require invertibility of $\D$. However, once we have 
proved Theorems 5.1 and \ref{exactly} in subsection 5.6, we replace 
(in the last subsection 5.7) $(\A,\HH,\D)$ by its double 
to remove 
this hypothesis. 

The standing assumption for the rest of this Section is that $(\A,\HH,\D)$ is a 
$QC^\infty$ finitely summable spectral triple with dimension $q\geq 1$. We will 
at times assume isolated spectral dimension, but shall always be explicit
in those results that need this hypothesis. 
{\bf Until we reach subsection \ref{noninvtble} we shall also assume that 
$\D$ is boundedly invertible}.

\subsection{Preliminary identities}
We begin by recalling some basic identities for the expectations 
$\la\cdots\ra_{m,r,s}$.

\begin{lemma}\label{basicidentities} For $m\geq 0$ and operators $A_0,...,A_m$, 
$A_i\in OP^{k_i}$, with $k_0+\cdots+k_m-2m-1<2Re(r)$, we have for $1\leq j<m$: 
\bean &&-\la A_0,...,[\D^2,A_j],...,A_m\ra_{m,s,r}\nno
&=&\la A_0,...,A_{j-1}A_j,...,A_m\ra_{m-1,s,r}-\la 
A_0,...,A_jA_{j+1},...,A_m\ra_{m-1,s,r};\eean
while for $j=m$ we have:
\bean &&-\la A_0,...,A_{m-1},[\D^2,A_m]\ra_{m,s,r}\nno
&=&\la A_0,...,A_{m-1}A_m\ra_{m-1,s,r}-(-1)^{A deg(A_m)}\la 
A_mA_0,...,A_{m-1}\ra_{m-1,s,r};\eean
We also have for $k\geq 1$:
\be \int_0^{\infty}s^k\la\D A_0,A_1,...,A_m\ra_{m,s,r}ds=
(-1)^A \int_0^{\infty}s^k\la A_0,A_1,...,A_m\D\ra_{m,s,r}ds.
\label{traceofcommutator}\ee
Moreover these $s$-integrals have a cyclic property:
$$\int_0^{\infty}s^k\la A_0,...,A_m\ra_{m,s,r}ds=(-1)^{A deg(A_m)}
\int_0^\infty s^k\la A_m,A_0,...,A_{m-1}\ra_{m,s,r}ds.$$
Furthermore if $\sum_{i=0}^m deg(A_i)\equiv A(\bmod\ 2)$, and we define 
$$deg_{-1}=0\;\;and\;\;deg_k=deg(A_0)+deg(A_1)+\cdots +deg(A_k)$$ then:
\be 
0=\sum_{j=0}^m(-1)^{deg_{j-1}}\int_0^{\infty}s^k\la A_0,...,[\D,A_j]_\pm,...,
A_m\ra_{m,s,r}ds.\label{addstozero}\ee
\end{lemma}

\begin{proof}
The first half of the first statement is just the first statement of 
Lemma 7.8 of \cite{CPRS2}
in the odd case: the proof in the even case is the same. The second part
of the first statement has a similar proof to the first part, but
requires the factor $(-1)^{Adeg(A_m)}$ since $A_m$ may anticommute with
$\gamma$ in the even case. The second statement
is also the second statement of Lemma 7.8 of \cite{CPRS2} in the odd case 
where $A=0.$
The even case is the same argument but needs the factor $(-1)^A$ since 
$\D$ and $\gamma$ anticommute. Similar comments apply to
the cyclic property, noting that $R_s(\lambda)$ has grading degree $0$
in either case. The cyclic property in the odd case is Lemma 7.7 of
 \cite{CPRS2}.\\
\hspace*{.2 in}The sum in (4) telescopes since 
$[\D,A_j]_\pm =\D A_j-(-1)^{deg(A_j)} A_j\D$ and $deg_j=deg_{j-1}+deg(A_j).$
After telescoping the sum, one is left with equation (3) since
$\sum_{i=0}^m deg(A_i)\equiv A(\bmod\ 2).$
\end{proof}

To introduce the `transgression cochain' we need a new expectation, and some 
basic properties.

\begin{defn} \label{newexpectation} 
For $m\geq 0$ and operators $A_0,...,A_m$, $A_i\in OP^{k_i}$, with 
$k_0+\cdots+k_m-2m-1<2Re(r)$  
define
\bean &&\la\la A_0,...,A_m\ra\ra_{m,s,r}\nno
&=&\sum_{j=0}^m(-1)^{deg_j}\la A_0,...,A_j,\D,A_{j+1},...,A_m\ra_{m+1,s,r}\nno
&=&\sum_{j=0}^m(-1)^{deg_j}\tau\left(\frac{1}{2\pi i}\; \gamma\int_l 
\lambda^{-q/2-r}
A_0R_s(\lambda)A_1\cdots A_jR_s(\lambda)\D R_s(\lambda)\cdots 
A_mR_s(\lambda)d\lambda\right).\eean
We note that except for a factor of $2$ and the possible $\pm 1$ 
factors this is just formal differentiation with respect to the ``variable''
$\D.$
\end{defn}

\begin{lemma} \label{lala-la-identity}For $m\geq 0$ and 
operators $A_0,...,A_m$, $A_i\in OP^{k_i}$, 
with $k_0+\cdots+k_m-2m-2<2Re(r)$ we have for $1\leq j < m,$ the identity:
\bea &&-\la\la A_0,...,[\D^2,A_j],...,A_m\ra\ra_{m,s,r}\nno
 &=&\la\la A_0,...,A_{j-1}A_j,...,A_m\ra\ra_{m-1,s,r}-\la\la 
A_0,...,A_jA_{j+1},...,A_m\ra\ra_{m-1,s,r}\nno
&+&(-1)^{deg_{j-1}}\la 
A_0,...,[\D,A_{j}]_\pm,...,A_m\ra_{m,s,r},\label{useforb}\eea
where we have a graded commutator in the last term. For $j=m$ we have the 
identity:
\bean &&-\la\la A_0,...,A_{m-1},[\D^2,A_m]\ra\ra_{m,s,r}\nno
 &=&\la\la A_0,...,A_{m-1}A_m\ra\ra_{m-1,s,r}-(-1)^{Pdeg(A_m)}\la\la 
A_mA_0,...,A_{m-1}\ra\ra_{m-1,s,r}\nno
&+&(-1)^{deg_{m-1}}\la 
A_0,...,[\D,A_{m}]_\pm\ra_{m,s,r}.\eean
If $\sum_{i=0}^mdeg(A_i)\equiv P(\bmod\ 2)$ and $\alpha\geq 1$ 
 we also have the  identity:
 \be \sum_{k=0}^m(-1)^{deg_{k-1}}\!\int_0^{\infty}\!\!s^\alpha\la\la 
A_0,..,[\D,A_k]_\pm,..,A_m\ra\ra_{m,s,r}ds=
\sum_{i=0}^m2\!\int_0^{\infty}\!\!s^\alpha\la 
A_0,..,A_i,\D^2,..,A_m\ra_{m+1,s,r}ds.\label{nothertrick}\ee
On the other hand, if $\sum_{i=0}^mdeg(A_i)\equiv A(\bmod\ 2)$ and $\alpha\geq 1$
then we have the following cyclic property for $\la\la\cdots\ra\ra$:
$$\int_0^\infty s^\alpha\la\la A_0,...,A_m\ra\ra_{m,s,r}ds=(-1)^{Pdeg(A_m)}
\int_0^\infty s^\alpha\la\la A_m,A_0,...,A_{m-1}\ra\ra_{m,s,r}ds.$$ 
><DEFANGED.836 \end{lemma}

\begin{proof} All statements are computations. The first two are careful 
applications of Lemma \ref{basicidentities} and are easily checked by 
the reader for say, $m=2$. The third is as follows 
(we suppress the integrals $\int_0^{\infty}s^\alpha...ds$ in order to 
lighten the
notation):
\bean && \sum_{k=0}^m(-1)^{deg_{k-1}}
\la\la A_0,...,[\D,A_k]_\pm,...,A_m\ra\ra_{m,s,r}\nno
&=&\sum_{k=0}^m(-1)^{deg_{k-1}}\sum_{i=0}^m\left\{
\begin{array}{ll}(-1)^{deg_i}\la 
A_0,...,A_i,\D,...,[\D,A_k]_\pm,...,A_m\ra_{m+1,s,r} 
& i<k\\ 
(-1)^{deg_i+1}\la A_0,...,[\D,A_k]_\pm,...,A_i,\D,...,A_m\ra_{m+1,s,r} & i\geq k
\end{array}\right.\nno
&=&\sum_{i=0}^m(-1)^{deg_i+1}\sum_{k=0}^m\left\{
\begin{array}{ll} (-1)^{deg_{k-1}+1}
\la A_0,...,A_i,\D,...,[\D,A_k]_\pm,...,A_m\ra_{m+1,s,r} & i<k\\ 
(-1)^{deg_{k-1}}\la A_0,...,[\D,A_k]_\pm,...,A_i,\D,...,A_m\ra_{m+1,s,r} & 
i\geq 
k
\end{array}\right.\nno
&+&\sum_{i=0}^m(-1)^{deg_i+1}\left((-1)^{deg_i}
\la A_0,...,A_i,2\D^2,...,A_m\ra_{m+1,s,r}\right.\nno
&&\qquad\qquad\qquad\qquad\qquad\qquad\qquad\qquad\qquad-\left.(-1)^{deg_i}\la 
A_0,...,A_i,2\D^2,...,A_m\ra_{m+1,s,r}\right)\nno
&=&\sum_{i=0}^m2\la A_0,...,A_i,\D^2,...,A_m\ra_{m+1,s,r}.\eean
The last line follows by applying Equation (\ref{addstozero}) above once for 
each $i=0,...,m$ and noting that 
the graded commutator of $\D$ with itself is $2\D^2$ and also that
$$\sum_{j=0}^m deg(A_j)  +  deg(\D) \equiv A (\bmod\ 2) \;\Longleftrightarrow\; 
\sum_{j=0}^m deg(A_j) \equiv P (\bmod\ 2).$$
The fourth identity is a calculation using the cyclic property for 
$\la\cdots\ra$ in Lemma 5.3.
\end{proof}
The next identity we refer to as the $\boldsymbol s$-{\bf trick}.

\begin{lemma} \label{strick} For any integers $m\geq 0, \alpha\geq 1$
and 
operators $A_0,...,A_m$ with $A_j\in OP^{k_j}$, $1+\alpha+\sum k_j-2m<2Re(r)$,
\ben \alpha\int_0^\infty s^{\alpha-1}\la 
A_0,...,A_m\ra_{m,s,r}ds=-2\sum_{j=0}^m\int_0^\infty s^{\alpha+1}\la 
A_0,...,A_j, 1, A_{j+1},...,A_m\ra_{m+1,s,r}ds.\een
\end{lemma}

{\bf Remark}. Provided $2Re(r)>\alpha+\sum k_j-2m$, the $s$-trick works exactly 
the same for the expectation $\la\la\cdots\ra\ra$.

\begin{proof} The derivative of $s^\alpha\la A_0,...,A_m\ra_{m,s,r}$ is 
\ben \alpha s^{\alpha-1}\la A_0,...,A_m\ra_{m,s,r}+2s^{\alpha+1}\sum_{i=0}^m\la 
A_0,...,A_i,1,A_{i+1},...,A_m\ra_{m+1,s,r}.\een 
Integrating this total derivative in $s$ from $0$ to $\infty$ gives the result.
\end{proof}

Similarly, we may employ our other variable of integration to obtain a 
$\boldsymbol\lambda$-{\bf trick}. The proof is virtually identical to the 
$s$-trick.

\begin{lemma}\label{lambdatrick} For any integer $m\geq0$, 
operators $A_0,...,A_m$ with $A_j\in OP^{k_j}$, and $r$ such that $\sum 
k_j-2m<2Re(r)$, we have
\ben -(q/2+r)\la A_0,...,A_m\ra_{m,s,r+1}=\sum_{j=0}^m\la 
A_0,...,A_j, 1, A_{j+1},...,A_m\ra_{m+1,s,r}.\een
\end{lemma}

\begin{proof}
We compute:
\begin{eqnarray*}
\lefteqn{\frac{d}{d\lambda}\left(\lambda^{-(q/2+r)}A_0R_s(\lambda)\cdots 
A_mR_s(\lambda)\right)}\\
&=& -(q/2+r)\lambda^{-(q/2+r+1)}A_0R_s(\lambda)\cdots A_mR_s(\lambda)\\
&&-\lambda^{-(q/2+r)}\sum_{j=0}^m A_0R_s(\lambda)\cdots A_j
R_s(\lambda)^2A_{j+1}
\cdots A_mR_s(\lambda).
\end{eqnarray*}
Integrating this equation along the line $l$, multiplying by 
$\frac{\gamma}{2\pi i}$ and taking the trace $\tau$ gives the result.
\end{proof}

\begin{lemma}\label{differentfort} For $m,\alpha\geq 0$, operators $A_i\in 
OP^{k_i}$ and $r$ such that $2Re(r)>1+\alpha-2m+\sum k_j$ we have
\bean &&\sum_{j=0}^m\int_0^\infty s^\alpha\la 
A_0,...,A_j,\D^2,A_{j+1},...,A_m\ra_{m+1,s,r}ds\nno
&=&-(m+1)\int_0^\infty s^\alpha\la 
A_0,...,A_m\ra_{m,s,r}ds+(1-q/2-r)\int_0^\infty s^\alpha\la 
A_0,...,A_m\ra_{m,s,r}ds\nno
&&+\frac{(\alpha+1)}{2}\int_0^\infty s^\alpha\la A_0,...,A_m\ra_{m,s,r}ds\nno
&&-\sum_{j=0}^m\int_0^\infty s^\alpha\la 
A_0,...,A_j,1,A_{j+1},...,A_m\ra_{m+1,s,r}ds\eean
\end{lemma}

\begin{proof}
This uses $\D^2R_s(\lambda)=-1+(\lambda-(1+s^2))R_s(\lambda)$. So
\bean && \la A_0,...,A_j,\D^2,A_{j+1},...,A_m\ra_{m+1,s,r}\nno
&=&-\la A_0,...,A_m\ra_{m,s,r}+\la 
A_0,...,A_j,1,A_{j+1},...,A_m\ra_{m+1,s,r-1}\nno
&&-(1+s^2)\la A_0,...,A_j,1,A_{j+1},...,A_m\ra_{m+1,s,r}.\eean
Summing gives us
\bean &&\sum_{j=0}^m\la A_0,...,A_j,\D^2,A_{j+1},...,A_m\ra_{m+1,s,r}\nno
&=&-(m+1)\la A_0,...,A_m\ra_{m,r,s}+(1-q/2-r)\la A_0,...,A_m\ra_{m,s,r}\nno
&&-\sum_{j=0}^m(1+s^2)\la A_0,...,A_j,1,A_{j+1},...,A_m\ra_{m+1,s,r}.\eean
here we used the $\lambda$-trick, Lemma 5.7.
Now integrate over $s$ and use the $s$-trick (Lemma 5.6) on the $s^2$ part 
of the last term to obtain the result.
\end{proof}

\subsection{The first homotopy}

We introduce a homotopy parameter 
into our resolvent cocycle and prove a transgression type formula.

\begin{defn} Assume that $\D$ is invertible. For $t\in[0,1]$ let 
\ben R_{s,t}(\lambda)=(\lambda-(t+s^2+\D^2))^{-1}\een
and
\ben \la\cdots\ra_{m,r,s,t},\ \ \ \la\la\cdots\ra\ra_{m,r,s,t}\een
be the expectations 
of Definitions \ref{expectation} and  \ref{newexpectation} using 
$R_{s,t}(\lambda)$ instead of $R_s(\lambda)$.
\end{defn}

To see that this is well-defined for $t\in[0,1]$, we must check that the trace 
estimates we require are still satisfied. These estimates all rest on the 
scalar inequality (see the proof of \cite[Lemma 5.2]{CPRS2}
$$(X+Y)^{-a-b}\leq X^{-a}Y^{-b}$$
for positive real numbers $X,Y,a,b.$
For invertible $\D$, there exists a $\delta>0$ such that 
$\D^2-\delta\geq 0$. Therefore $\D^{-2}$, $(\D^2-\delta/2)^{-1}$, 
and $(1+\D^2)^{-1}$ all have the same summability (say, for any $K>q/2$).
So by the functional calculus
if $K\geq q/2 + \epsilon$ we have for all $t\in [0,1]$:
\bean 
(\D^2+t+s^2)^{-K}&=&(\D^2-\delta/2+\delta/2+t+s^2)^{-K}
\leq(\D^2-\delta/2)^{-q/2-\epsilon}
(\delta/2+t+s^2)^{-K+q/2+\epsilon}\nno
&\leq &(\D^2-\delta/2)^{-q/2-\epsilon}
(\delta/2+s^2)^{-K+q/2+\epsilon}.\eean
This is enough to prove part $(b)$ of the following lemma. Part $(a)$
is just the functional calculus. It is important to note that these lemmas
imply that all the trace estimates will continue to hold uniformly 
for $t\in[0,1]$ provided that $\D$ is invertible. For our $\lambda$-integral we 
choose our line $l$ given by $\lambda=a+iv$ so that $0<a<\delta/4$ to
simply certain estimates, although this is not always needed.

In particular, we will need the following modifications of the results in
\cite{CPRS2}.

\begin{lemma}\label{normlambda} (Cf \cite[Lemma 5.1]{CPRS2}) Let $\D$
be an unbounded self-adjoint operator with $\D^2\geq\delta$ .\\
(a) For $\lambda=a+iv\in \C$, $0<a<\delta/4$, and $t,s\geq 0$ we have the 
estimate 
\ben \n(\lambda-(t+\D^2+s^2))^{-1}\n\leq (v^2+(s^2+t+\delta -a)^2)^{-1/2}
\leq\frac{1}{\delta -a}.\een
(b) For $t,s\geq 0$ we have the 
estimate 
\ben (\D^2+t+s^2)^{-K}
\leq (\D^2-\delta/2)^{-q/2-\epsilon}
(\delta/2+s^2)^{-K+q/2+\epsilon}.\een
\end{lemma}

\begin{lemma}\label{lambda} (Cf \cite[Lemma 5.3]{CPRS2}) Let $\D^2\geq\delta.$
Let $\lambda=a+iv\in \C$, $0<a<\delta/4$, and $t,s\geq 0.$ 
For $q\geq 1$ let $(1+\D^2)^{-1/2}$
be $(q+\epsilon)$-summable for every $\epsilon >0$. 
Then for each $\epsilon >0$ and 
$N>(q+\epsilon)/2$ , we have the trace-norm estimate:  
\begin{eqnarray*}\n (\lambda-(t+\D^2+s^2))^{-N}\n_1&\leq& C'_{q+\epsilon}
((\delta/2+s^2-a)^2+v^2)^{-N/2+(q+\epsilon)/4}\\
&\leq&C'_{q+\epsilon}((s^2+a)^2+v^2)^{-N/2+(q+\epsilon)/4}.
\end{eqnarray*}
where $C'_{q+\epsilon}= \n (\D^2-\delta/2)^{-(q+\epsilon)/2}\n_1.$
\end{lemma}

\begin{lemma}\label{intest}(Cf \cite[Lemma 5.4]{CPRS2}
 Let $0<a<\delta/4$, $t\in [0,1]$ and let $J,K,B\geq 0$. 
Then 
\begin{eqnarray*}
&&\int_0^\infty \int_{-\infty}^\infty s^J\sqrt{a^2+v^2}^{-B}
\sqrt{(t+s^2+\delta/2-a)^2+v^2}^{-K}dvds\\
&\leq&\int_0^\infty \int_{-\infty}^\infty s^J\sqrt{a^2+v^2}^{-B}
\sqrt{(s^2+a)^2+v^2}^{-K}dvds
\end{eqnarray*}
converges provided $J-2K<-1$ and $J-2K+1-2B<-2$. Moreover, if $K$ and $J$
are fixed with $J-2K<-1$, then the integral as a function of $B$ goes to
$0$ as $B\to\infty$. In particular, with $J,K$ fixed and $B_0$ some value 
for which the integral is finite, then the integrals are uniformly bounded
for all $B\geq B_0.$
\end{lemma}

\begin{defn}\label{transgressioncochain} Define, for $Re(r)>(1-m)/2$, 
the components of the {\bf `transgression'} 
cochain
\ben \Phi^r_{m,t}(a_0,...,a_m)=\frac{\eta_{m+1}}{2}\int_0^\infty s^{m+1}\la\la 
a_0,da_1,...,da_m\ra\ra_{m,s,r,t}ds.\een
Similarly, we define $\phi^r_{m,t}$ using $R_{s,t}(\lambda)$ in place of
$R_s(\lambda)$ in the definition of $\phi^r.$
\end{defn}

\begin{rems*}
Lemmas 5.3 to 5.7 inclusive work equally well for the expectations 
$\la\cdots\ra_{m,s,r,t}$ and $\la\la\cdots\ra\ra_{m,s,r,t}$, however Lemma 5.8
has an extra factor of $t$ on the last term as we will make explicit at the end
of the proof of the following proposition.
\end{rems*}

\begin{prop}\label{tcohom} We have the $(b,B)$ bicomplex formula (with 
coefficients in the functions holomorphic for $Re(r)>1/2$) for 
$m\equiv P(\bmod\ 2)$
\ben 
(B\Phi^r_{m+1,t}+b\Phi_{m-1,t}^r)(a_0,...,a_m)=\left(\frac{q-1}{2}+r\right)
\phi_{m,t}^r(a_0,...,a_m)-t\frac{q+2r}{2}\phi^{r+1}_{m,t}(a_0,...,a_m).\een
\end{prop}

\begin{proof}
We begin the proof by computing with $\Phi^r_{m}:=\Phi^r_{m,1}$. First,
using  the cyclic property of $\la\la\cdots\ra\ra$ of Lemma 5.5 and the fact 
that $m\equiv P(\bmod\ 2),$ we have writing $da=[\D,a]$:
\bea B\Phi^r_{m+1}(a_0,...,a_m)&=&\frac{\eta_{m+2}}{2}\sum_{j=0}^m\int_0^\infty 
s^{m+2}(-1)^{mj}\la\la 1,da_j,...,da_{j-1}\ra\ra_{m+1,s,r} ds\nno
&=&\frac{\eta_{m+2}}{2}\sum_{j=0}^m\int_0^\infty s^{m+2}\la\la 
da_0,...,da_{j-1},1,da_j,...,da_m\ra\ra_{m+1,s,r} ds\nno
&=&-\frac{\eta_{m+2}(m+1)}{4}\int_0^\infty s^{m}\la\la 
da_0,...,da_m\ra\ra_{m,s,r} ds\nno
&=&-\frac{\eta_{m}}{2}\int_0^\infty s^{m}\la\la 
da_0,...,da_m\ra\ra_{m,s,r} ds\label{formulaforBPhi}\eea
using the $s$-trick (Lemma 5.6) in the last line, which is the same for 
$\la\la\cdots\ra\ra$ as for $\la\cdots\ra$.

The computation for $b\Phi_{m-1}^r$ is the same as for $b\phi^r_{m-1}$, 
\cite{CPRS2}, except we need to take account of Equation (\ref{useforb}). This 
gives
\bean 
b\Phi_{m-1}^r(a_0,...,a_m)&=&\frac{\eta_{m}}{2}\sum_{j=1}^m(-1)^j\int_0^\infty 
s^m\la\la a_0,da_1,...,[\D^2,a_j],...,da_m\ra\ra_{m,s,r} ds\nno
&&-\frac{\eta_{m}}{2}\sum_{j=1}^m\int_0^\infty s^m\la 
a_0,da_1,...,da_j,...,da_m\ra_{m,s,r} ds\nno
&=&\frac{\eta_{m}}{2}\sum_{j=1}^m(-1)^j\int_0^\infty s^m\la\la 
a_0,da_1...,[\D^2,a_j],...,da_m\ra\ra_{m,s,r} ds\nno
&&-\frac{\eta_{m}m}{2}\int_0^\infty s^m\la 
a_0,da_1,...,da_m\ra_{m,s,r} ds.\eean

Now put them together. First, using $\eta_{m+2}(m+1)/2=\eta_{m}$ we get
\bean &&(B\Phi^r_{m+1}+b\Phi_{m-1}^r)(a_0,...,a_m)\nno
&=&-\frac{\eta_m}{2}\int_0^\infty s^{m}\la\la 
da_0,...,da_m\ra\ra_{m,s,r} ds\nno
&+&\frac{\eta_{m}}{2}\sum_{j=1}^m(-1)^j\int_0^\infty s^m\la\la 
a_0,da_1,...,[\D^2,a_j],...,da_m\ra\ra_{m,s,r} ds\nno
&&-\frac{\eta_{m}m}{2}\int_0^\infty s^m\la 
a_0,da_1,...,da_m\ra_{m,s,r} ds\eean

then using the fact that $[\D^2,a_j]=[\D,[\D,a_j]]_{\pm}$ we get:

\bean
&=&-\frac{\eta_m}{2}(-1)^{deg(a_0)}\int_0^\infty s^{m}\la\la 
[\D,a_0]_\pm,da_1,...,da_m\ra\ra_{m,s,r} ds\nno
&+&\frac{-\eta_{m}}{2}\sum_{j=1}^m(-1)^{deg(a_0)+deg(da_1)+\cdots+deg(da_
{j-1})}\int_0^\infty s^m\la\la a_0,da_1...,[\D,da_j]_\pm,...,da_m\ra\ra_{m,s,r} 
ds\nno
&&-\frac{\eta_{m}m}{2}\int_0^\infty s^m\la 
a_0,da_1,...,da_m\ra_{m,s,r} ds\eean

now, by identity (\ref{nothertrick}) of Lemma 5.5, this gives:
\bean
&=&\frac{-2\eta_{m}}{2}\sum_{j=0}^m\int_0^\infty s^m\la 
a_0,...,da_j,\D^2,da_{j+1},...,da_m\ra_{m+1,s,r} ds\nno
&&-\frac{\eta_{m}m}{2}\int_0^\infty s^m\la 
a_0,da_1,...,da_m\ra_{m,s,r} ds\eean

then, applying Lemma 5.8 we get:
\bea &=&-\eta_{m}\left((-(m+1)\int_0^\infty s^m\la 
a_0,da_1,...,da_m\ra_{m,s,r}ds+\right.\nno
&+&\left.(1-q/2-r)\int_0^\infty s^m\la 
a_0,da_1,...,da\ra_{m,s,r}ds\right.\nno
&&+\frac{(m+1)}{2}\int_0^\infty s^m\la 
a_0,da_1,...,da_m\ra_{m,s,r}ds\nno
&&\left.-\sum_{j=0}^m\int_0^\infty s^m\la 
a_0,...,da_j,1,da_{j+1},...,da_m\ra_{m+1,s,r}ds\right)\nno
&&-\frac{\eta_{m}m}{2}\int_0^\infty s^m\la 
a_0,da_1,...,da_m\ra_{m,s,r} ds\nno
&=&\eta_{m}\frac{q+2r-1}{2}\int_0^\infty s^m\la 
a_0,da_1,...,da_m\ra_{m,s,r} ds\nno
&&+\eta_{m}\sum_{j=0}^m\int_0^\infty s^m\la 
a_0,...,da_j,1,da_{j+1},...,da_m\ra_{m+1,s,r}ds\nno
&=&\frac{q+2r-1}{2}\phi^r_m(a_0,...,a_m)-\frac{q+2r}{2}\phi_m^{r+1}(a_0,...,a_m)
.\label{aaagh}\eea
We used the $\lambda$-trick (Lemma 5.7) in the last line. 

We now do the general case of $t\in [0,1]$ and observe that for $\alpha\geq 1$, 
a slight variation on Lemma \ref{differentfort} gives
\bean &&\sum_{j=0}^m\int_0^\infty s^\alpha\la 
A_0,...,A_j,\D^2,A_{j+1},...,A_m\ra_{m+1,r,s,t}ds\nno
&=&-(m+1)\int_0^\infty s^\alpha\la A_0,...,A_m\ra_{m,r,s,t}ds+(1-q/2-r)
\int_0^\infty s^\alpha\la A_0,...,A_m\ra_{m,r,s,t}ds\nno
&&+\frac{(\alpha +1)}{2}\int_0^\infty s^\alpha\la A_0,...,A_m\ra_{m,r,s,t}ds\nno
&&-t\sum_{j=0}^m\int_0^\infty s^\alpha\la 
A_0,...,A_j,1,A_{j+1},...,A_m\ra_{m+1,r,s,t}ds.\eean
As mentioned in the remark, the use of Lemmas 5.3 to 5.7 in the proof for $t=1$ 
work equally well for $t\neq 1$ and now Lemma 5.8 is modifed as above. Putting
these results together with the computation in the case $t=1$ gives us:

\bean 
(B\Phi^r_{m+1,t}+b\Phi_{m-1,t}^r)(a_0,...,a_m)&=&
\left(\frac{q+2r-1}{2}\right)\phi_{m,t}^r(a_0,...,a_m)\nno
&+& t \eta_{m}\sum_{j=0}^m\int_0^\infty s^m\la 
a_0,...,da_j,1,da_{j+1},...,da_m\ra_{m+1,r,s,t}ds\nno
&=&\left(\frac{q+2r-1}{2}\right)\phi_{m,t}^r(a_0,...,a_m)-t\frac{q+2r}{2}
\phi^{r+1}_{m,t}(a_0,...,a_m).\eean
\end{proof}

\begin{cor}\label{tequalszero} For $\D$ invertible and $m\equiv
P(\bmod\ 2)$ (setting $t=0$)
 we have
\ben 
(B\Phi^r_{m+1,0}+b\Phi_{m-1,0}^r)(a_0,...,a_m)=\left(\frac{q-1}{2}+r\right)\phi_
{m,0}^r(a_0,...,a_m)\een
\end{cor}

\begin{prop}\label{singleterm} For $\D$ invertible, the cochain
$(\phi^r_{m,0})_{m=P}^\infty$ is 
cohomologous to the  cochain,
\be \frac{1}{(r-(1-q)/2)}B\Phi^r_{M+1,0}\label{almostchern}.\ee
This 
cochain (\ref{almostchern}) is a cyclic cocycle modulo cochains with values in 
the functions holomorphic at $r=(1-q)/2$.
\end{prop}

\begin{proof} Let $M=2N-P$. By Proposition \ref{tequalszero} the cochain 
given by the infinite tuple
\ben 
\left(\frac{1}{(r+(q-1)/2)}\Phi_{P-1,0}^r,\dots,\frac{1}{(r+(q-1)/2)}
\Phi^r_{M-1,0},
0,0,\dots\right)\een
$(B,b)$-cobounds the following difference. That is, applying $(B,b)$ to the 
above cochain yields:
\bean 
&&\left(\phi^r_{P,0},\phi^r_{P+2,0},\dots ,\phi^r_{M,0}-\frac{B\Phi^r_{M+1,0}}
{(r+(1-q)/2)},0,0,\dots\right)=\left((\phi^r_{m,0})_{m=P}^{M}-
\frac{B\Phi^r_{M+1,0}}{(r+(q-1)/2)}\right).
\eean
That is, $(\phi^r_{m,0})$ is cohomologous to 
$\frac{1}{(r+(q-1)/2)}B\Phi^r_{M+1,0}$. Observe that 
 because it is in the image of $B$, 
$(r-(1-q)/2)^{-1}B\Phi^{r}_{M+1,0}$ is cyclic. It is also a $b$-cyclic cocycle 
modulo cochains with values in the functions holomorphic at $r=(1-q)/2$. This 
follows from 
\ben b\Phi^r_{M-1,0}+B\Phi^r_{M+1,0}=(r+(q-1)/2)\phi^r_{M,0}\een
by applying $b$ and recalling that $b\phi^r_{M,0}$ is holomorphic at 
$r=(1-q)/2$, \cite{CPRS2}. 
\end{proof}

\begin{cor}\label{isolatedcohom} For $\D$ invertible and assuming that 
$(\A,\HH,\D)$ has isolated spectral dimension, the residue cocycle is 
cohomologous to $B\Phi^{(1-q)/2}_{M+1,0}$. Moreover $B\Phi^{(1-q)/2}_{M+1,0}$  
is a cyclic cocycle.
\end{cor}

\begin{proof}
By Proposition \ref{tequalszero} we have  
\ben 
\left(res_{r=(1-q)/2}\frac{1}{(r+(q-1)/2)}\Phi_{P-1,0}^r,\dots,res_{r=(1-q)/2}
\frac{1}{(r+(q-1)/2)}\Phi^r_{M-1,0},0,0,\dots\right)\een
$(B,b)$-cobounds the difference
\bean 
&&res_{r=(1-q)/2}\left(\phi^r_{P,0},\phi^r_{P+2,0},\dots ,\phi^r_{M,0}-
\frac{B\Phi^r_{M+1,0}}{(r+(1-q)/2)},0,0,\dots\right)\nno
&=&\left(\phi_{P,0},\phi_{P+2,0},\dots ,\phi_{M,0}-B\Phi^{(1-q)/2}_{M+1,0},0,0,
\dots\right).\eean
This statement requires that $b\Phi^r_{M-1,0}(a_0,...,a_M)$ is holomorphic at 
$r=(1-q)/2$, and that $B\Phi_{M+1,0}^r$ is finite (and so holomorphic) at 
$r=(1-q)/2$. This is easy to prove using \cite[Lemmas 7.1,
7.2]{CPRS2}. 
Finally, to see that 
$B\Phi^{(1-q)/2}_{M+1,0}$ is a cyclic cocycle, we simply take residues of the 
corresponding result in Proposition \ref{singleterm}.
\end{proof}

\subsection{The homotopy to $t=0$.}

In what follows we suppose that $t,t+\epsilon\in[0,1]$, and we write
$$ R_{s,t}(\lambda)=(\lambda-(t+s^2+\D^2))^{-1}.$$
By the resolvent equation: $R_{s,t+\epsilon}(\lambda)-R_{s,t}(\lambda)
=\epsilon R_{s,t+\epsilon}(\lambda)R_{s,t}(\lambda)$ and the fact that
$|R_{s,t}(\lambda)|\leq |R_{s,0}(\lambda)|$ we see that
$$\n R_{s,t+\epsilon}(\lambda)-R_{s,t}(\lambda)\n_{M+1}\leq
\frac{\epsilon}{\delta -a}\n R_{s,0}(\lambda)\n_{M+1}.$$
Similar considerations show that:
$$\n\frac{1}{\epsilon}(R_{s,t+\epsilon}(\lambda)-R_{s,t}(\lambda))
-R_{s,t}(\lambda)^2\n_{M+1}\leq
\frac{\epsilon}{(\delta-a)^2}\n R_{s,0}(\lambda)\n_{M+1}.$$
At this point we lighten the notation temporarily by dropping the
dependence on $\lambda.$ More generally if $n\geq 1$, we can use the identity
$X^n-Y^n=\sum_{k=0}^{n-1}X^k(X-Y)Y^{n-k-1}$ and H\"{o}lder's inequality
to show that:
$$\n R_{s,t+\epsilon}^n-R_{s,t}^n\n_{\frac{M+1}{n}}\leq\frac{\epsilon\cdot n}
{\delta-a}\n R_{s,0}\n_{M+1}^n\;\;\;and$$
$$\n\frac{1}{\epsilon}(R_{s,t+\epsilon}^n-R_{s,t}^n) -nR_{s,t}^{n+1}
\n_{\frac{M+1}{n}}\leq\frac{\epsilon}{(\delta-a)^2}\cdot\frac{n(n+1)}{2}
\n R_{s,0}\n_{M+1}^n.$$

Let $F_m$ be the space of functions defined and holomorphic in the
right half plane $\{z\in\C:Re(z)>(1-m)/2\}$, and give $F_m$ the
topology of uniform convergence on compacta.

\begin{prop}\label{generalt} For each $m=P,P+2,\dots$, the map
$$[0,1]\ni t\mapsto \phi^\bullet_{m,t}\in H(\A^{\otimes m+1},F_m)$$
is $C^1$ and 
$$\frac{d}{dt}\phi^\bullet_{m,t}=-(q/2+\bullet)\phi^{\bullet+1}_{m,t}.$$
\end{prop}

\begin{proof} We do the case $m<M$ where we must use some initial trickery 
to get to a computable situation. For $m\geq M$ such tricks are not needed.
We recall from the proof of Lemma \ref{better-be-true} that if we apply the
$s$-trick $(M-m)/2$ times and the $\lambda$-trick $(M-m)/2$ times we get
by summing over all $k=(k_0,\dots,k_m)$ where each $0\leq k_i\leq (M-m)$ 
and $|k|=(M-m)$:
\begin{eqnarray*}
&&\int_0^\infty s^m \la a_0,da_1,\dots,da_m\ra_{m,s,r,t} ds
=(2)^{\frac{M-m}{2}}(M-m)!\prod_{l=1}^{\frac{M-m}{2}}\frac{1}{q/2+r-l}
\prod_{j=1}^{\frac{M-m}{2}}\frac{1}{m+j}\times\\
&&\times\sum_k\int_0^\infty s^M
\la a_0,1^{k_0},da_1,1^{k_1},\dots,da_m,1^{k_m}\ra_{M,s,r-(M-m)/2,t} ds
\end{eqnarray*}
Where, we mean $1^{k_i}=1,1,\dots,1$ with $k_i$ one's. Ignoring the prefactors
on the right hand side, this becomes:
$$\sum_n\int_0^\infty s^M\tau\left(\frac{1}{2\pi i}\gamma
\int_l\lambda^{-q/2-r+(M-m)/2}
a_0R_{s,t}^{n_0}da_1R_{s,t}^{n_1}\cdots da_mR_{s,t}^{n_m} d\lambda\right)ds$$
where we sum over $n=(n_0,\dots,n_m)$ with each $n_i=k_i+1$
so that $1\leq n_i\leq(M-m)+1$ and $|n|=M+1.$ Now each integrand is not only 
trace-class, but by the estimates immediately preceding the statement of this 
proposition they are $t$-differentiable in trace-norm using the usual product 
rule argument and H\"{o}lder's inequality. In particular we have:
$$\n\frac{1}{\epsilon}\left(a_0R_{s,t+\epsilon}^{n_0}da_1\cdots 
da_mR_{s,t+\epsilon}^
{n_m}-a_0R_{s,t}^{n_0}da_1\cdots da_mR_{s,t}^{n_m}\right)-
\sum_{i=0}^m n_ia_0R_{s,t}^{n_0}\cdots da_iR_{s,t}^{n_i+1}\cdots 
da_mR_{s,t}^{n_m}\n_1$$
$$\leq C\epsilon\n a_0\n\cdot\n da_1\n\cdots\n da_m\n\cdot\n R_{s,0}\n_{M+1}.$$
Where the constant $C$ is independent of $s$, $\lambda$ and $r$. Since
$$\int_0^\infty s^M\int_{l} |\lambda^{-q/2-r+(M-m)/2}|\cdot
\n R_{s,0}(\lambda)\n_{M+1} d\lambda ds = O(Re(r)^{-1}) < \infty$$
by Lemma \ref{intest}, we can invoke the Lebesgue Convergence Theorem to 
conclude that:\\ $d/dt(\phi_{m,t}^r(a_0,...,a_m))$ exists and equals
$$\eta_m 2^{\frac{M-m}{2}}(M-m)!\prod_{b=1}^{\frac{M-m)}{2}}\frac{1}{q/2+r-b}
\prod_{j=1}^{\frac{M-m)}{2}}\frac{1}{m+j}\times$$
$$\times \sum_k\sum_{i=0}^m\int_0^\infty s^M (k_i+1)\la a_0,1^{k_0},\dots,
da_i,1^{k_i+1},\dots,da_m,1^{k_m}\ra_{M+1,s,r-(M-m)/2,t} ds.$$
Now undoing our applications of the $s$-trick and the $\lambda$-trick we get:
$$ \frac{d}{dt}\phi^r_{m,t}(a_0,\dots,a_m)=
\eta_m\sum_{j=0}^m\int_0^\infty s^m\la a_0,\dots,da_j,1
,da_{j+1},\dots,da_m \ra_{m+1,s,r,t}ds$$
and a final application of the $\lambda$-trick yields:
$$ \frac{d}{dt}\phi^r_{m,t}(a_0,\dots,a_m)=
-(q/2+r)\phi^{r+1}_{m,t}(a_0,\dots,a_m).$$ 
We note that by our estimates the convergence is uniform in $r$ for $r$ in
a compact set.
\end{proof}

In the discussion below $k=0,1,2,\dots$, $m=P,P+2,\dots , M$, $t\in [0,1]$.

\begin{lemma}\label{iteration}
We have
$$\phi^{r+k}_{m,t}
=\frac{1}{r+k+(q-1)/2}[B\Phi^{r+k}_{m+1,t} + b\Phi^{r+k}_{m-1,t} 
+(\frac{q}{2} +r +k) t\phi^{r+k+1}_{m,t}]$$
\end{lemma}

\begin{proof} This is just Proposition \ref{tcohom} with $r+k$ in place of $r.$
\end{proof}

\begin{prop}\label{cohomforallt} For all $R,\ T\in[0,1]$, the cocycles 
$(\phi^r_{m,T})_{m=P}^M$ and $(\phi^r_{m,R})_{m=P}^M$
are equal modulo coboundaries and modulo cochains yielding functions 
holomorphic at $r=(1-q)/2$.
\end{prop}

\begin{proof}
Recall from Proposition \ref{generalt} that for $\D$ invertible,
$\phi^r_{m,t}$ is defined and holomorphic for $Re(r)>(1-m)/2$ for all
$t\in[0,1]$. As $[0,1]$ is compact, the integral
$$\int_0^1\phi^r_{m,t}(a_0,\dots,a_m)dt$$
is holomorphic for $Re(r)>(1-m)/2$ and any $a^0,\dots,a^m\in\A$. 

Now we make some simple observations, omitting the
variables $a_0,\dots,a_m$ to lighten the notation. For $T,R\in[0,1]$ we
have
\be\phi^r_{m,T}-\phi^r_{m,R}=\int_R^T \frac{d}{dt}\phi^r_{m,t}dt
=-(q/2+r)\int_R^T\phi^{r+1}_{m,t}dt.\label{dopeytrick}\ee

Now apply the formula of Lemma \ref{iteration} iteratively. First we get

$$\phi^r_{m,T}-\phi^r_{m,R}=\frac{-(q/2+r)}{r+1+(q-1)/2}
\int_R^T\left(B\Phi^{r+1}_{m+1,t} + b\Phi^{r+1}_{m-1,t} 
+(\frac{q}{2} +r +1) t\phi^{r+2}_{m,t}\right)dt.$$

Observe that the numerical factors are holomorphic at
$r=(1-q)/2$. 

Iterating this procedure $L$ times gives us
\bean\phi^r_{m,T}-\phi^r_{m,R}&=&\frac{-(q/2+r)\cdots(q/2+r+L)}
{(r+1+(q-1)/2)\cdots(r+L+(q-1)/2)}
\int_R^T t^L\phi^{r+L+1}_{m,t}dt\nno
&
+&\sum_{j=1}^L\frac{-(q/2+r)\cdots(q/2+r+j-1)}{(r+1+(q-1)/2)\cdots
(r+j+(q-1)/2)}\int_R^T\left(B\Phi^{r+j}_{m+1,t} + b\Phi^{r+j}_{m-1,t} 
\right)t^{j-1}dt.\eean
We would like to know (for completeness) what is the smallest integer $L$
that guarantees that $\phi^{r+L+1}_{m,t}$ is holomorphic at
$r=(1-q)/2$ for all $m$. We know that we require
$$Re(r)+L+1>(1-m)/2$$
for each $m$. Rearranging gives $Re(r)>-1-L+(1-m)/2$, and we would like
the right hand side to be strictly less than $(1-q)/2$ for each
$m$. In the even case the worst situation is when $m=0$, and for this
case we require
$L>-1+q/2$. Since in the even case $N=[(q+1)/2]>(q-1)/2>-1+q/2$, we
may choose $L=N$. In the odd case the worst case is $m=1$, and we
require $L>-1+(q-1)/2$. In this case $N=[(q+2)/2]>q/2>(q-1)/2$, and so
we may take $L=N-1$. (If $N-1=0$, then we need do no iterations,
because Equation \eqref{dopeytrick} shows that modulo cochains
yielding holomorphic functions, $\phi^r_{1,T}=\phi^r_{1,R}$, and these
are the only functionals in the resolvent cocycle).
 
With this choice of $L=N-P$, we have modulo cochains yielding functions 
holomorphic in a
half plane containing $(1-q)/2$,
$$\phi^r_{m,T}-\phi^r_{m,R}=\sum_{j=1}^L\frac{-(q/2+r)\cdots(q/2+r+j-1)}
{(r+1+(q-1)/2\cdots(r+j+(q-1)/2)}\int_R^T\left(B\Phi^{r+j}_{m+1,t}
  + b\Phi^{r+j}_{m-1,t} \right)t^{j-1}dt.$$

Thus a simple rearrangement yields the cohomology, valid for $Re(r)>(1-P)/2$,
\bean
&&(\phi^r_{m,T}-\phi^r_{m,R})_{m=P}^M-B\sum_{j=1}^L\frac{-(q/2+r)\cdots
(q/2+r+j-1)}{(r+1+(q-1)/2\cdots(r+j+(q-1)/2)}
\int_R^T\Phi^{r+j}_{M+1,t}t^{j-1}dt\nno
&=&
(B+b)\left(\sum_{j=1}^L\frac{-(q/2+r)\cdots(q/2+r+j-1)}{(r+1+(q-1)/2
\cdots(r+j+(q-1)/2)}\int_R^T\Phi^{r+j}_{m,t}t^{j-1}dt\right)_{m=P}^{M-1}.
\eean
Hence modulo coboundaries and cochains yielding functions holomorphic at
$r=(1-q)/2$, we have the equality
$$(\phi^r_{m,T}-\phi^r_{m,R})_{m=P}^M=B\sum_{j=1}^L\frac{-(q/2+r)
\cdots(q/2+r+j-1)}{(r+1+(q-1)/2\cdots(r+j+(q-1)/2)}\int_R^T
\Phi^{r+j}_{M+1,t}t^{j-1}dt.$$
However, an application of \cite[Lemma 7.2]{CPRS2} 
now shows that the right hand side is holomorphic at $r=(1-q)/2$,
since $j\geq 1$ in all cases. Hence,
modulo coboundaries and cochains yielding functions holomorphic at
$r=(1-q)/2$, we have the equality
$$(\phi^r_{m,T})_{m=P}^M=(\phi^r_{m,R})_{m=P}^M.$$
\end{proof}

\begin{cor}\label{wotwewanted} Modulo coboundaries and cochains 
yielding functions holomorphic in a half plane containing
$r=(1-q)/2$, we have the equality
$$(\phi^r_m)_{m=P}^M=:(\phi^r_{m,1})_{m=P}^M=B\Phi^r_{M+1,0}.$$
\end{cor}

\subsection{Homotopy to the Chern character for invertible $\D$}
We now drop the $0$ subscript from $\Phi^r_{M+1,0}$ as we will 
leave $t=0$ from 
now on, and consider a different homotopy. 
We now want to deform $B\Phi_{M+1}^{r}$ using the homotopy
$ u\tto \D_u:=\D\dd^{-u},$ following the strategy of Higson.
Unfortunately we discovered that this homotopy is quite subtle
if one wants to check all of the estimates needed to show it is well defined.
Hence the somewhat lengthy discussion in this subsection.
Of course it is clear that this homotopy only makes sense in the 
invertible case. 
We handle the transition to the non-invertible case later.

We write 
 $\Phi^r_{u,M+1}$ for $\Phi^r_{M+1}$ defined using 
$\D_u$ instead of $\D=\D_0$. We would also like to write $\dot\D=-\D_u\log\dd$, 
but this is purely formal and we will only take the limit of the corresponding 
difference quotients when they are multiplied by a factor of $|\D|^{-\rho}$
to ensure that the limit exists. 
In fact this derivative only appears via the derivative of $B\Phi^r_{u,M+1}$, 
which 
in turn appears in the coboundary computation of yet another auxiliary cochain 
$\Psi^r_{u,M}$ in Lemma \ref{aagghhh} below. The necessary estimates for taking 
the derivative of $B\Phi^r_{u,M+1}$ are proved in Lemma \ref{had-to-be-done}.

\begin{lemma}\label{aagghhh} For $r>(1-M)/2$ define a functional by 
\ben \Psi^r_{u,M}(a_0,...,a_M)=-\frac{\eta_M}{2}\int_0^\infty s^M\la\la 
a_0\dot{\D}_u,[\D_u,a_1],...,[\D_u,a_M]\ra\ra_{M,s,r} ds\een
Then we have
\bean &&(bB\Psi^r_{u,M})(a_0,...,a_M)\nno
&=&-\eta_M(r+(q-1)/2)\sum_{i=0}^M(-1)^i\int_0^\infty s^M\la 
[\D_u,a_0],...,[\D_u,a_i],\dot{\D}_u,...,[\D_u,a_M]\ra_{M+1,s,r} ds\nno
&&+\frac{d}{du}(B\Phi^r_{u,M+1})(a_0,...,a_M).\eean
\end{lemma}

{\bf Remarks}. (i) A judicious use of H\"older's inequality shows directly that
$\Psi^r_{u,M}$ is finite for $Re(r)>(1-M)/2$. However we choose a slightly 
different 
argument and obtain this fact as a corollary of  the computations
below               . 

(ii) The derivative of $B\Phi^r_{u,M+1}$ must be shown to exist for 
$Re(r)>(1-M)/2$. The functional $B\Phi^r_{u,M+1}$ is finite in this
region, 
and 
so a similar argument to that which shows $\phi^r_m$ is holomorphic when
finite  will 
show that  $B\Phi^r_{u,M+1}$ is 
holomorphic when finite. In fact $B\Phi^r_{u,M+1}$ is finite 
for $Re(r)>-M/2$ by \cite[Lemma 7.2]{CPRS2}. 
The problem is the
existence of the derivative in $u$, and for this we require 
a careful argument which we give in the next Lemma.

We will return to the proof of \ref{aagghhh} after we handle the
technical issues which are summarised in the next result.

\begin{lemma}\label{had-to-be-done} Write $R_u$ for the resolvent 
$R_u=(\lambda-(s^2+\D_u^2))^{-1}$ defined using $\D_u$. For 
$a_0,\dots,a_M\in\A$, the product
$$T_{u,j}:=[\D_u,a_0]R_u[\D_u,a_1]R_u\cdots
[\D_u,a_j]R_u\D_u R_u[\D_u,a_{j+1}]\cdots
 R_u[\D_u,a_M]R_u$$
is trace class for all $u\in[0,1]$ and is $C^1$ in $u$. 
Moreover, with $\dot{\D_u}=-\D_u log(|\D|)$
we get as expected:
\begin{eqnarray*}\frac{d}{du}(T_{u,j})
&=&\sum_{k=0}^M [\D_u,a_0]R_u\cdots R_u[\D_u,a_k]
(R_u 2\D_u\dot{\D_u}R_u)[\D_u,a_{k+1}]R_u\cdots[\D_u,a_M]R_u\\
&+& [\D_u,a_0]R_u\cdots R_u[\D_u,a_j]R_uD_u(R_u 2\D_u\dot{\D_u}R_u)[\D_u,a_{j+1}]
\cdots R_u[\D_u,a_M]R_u\\
&+& \sum_{k=0}^M [\D_u,a_0]R_u[\D_u,a_1]R_u\cdots R_u
[\dot{\D_u},a_k]R_u\cdots R_u[\D_u,a_M]R_u\\
&+& [\D_u,a_0]R_u[\D_u,a_1]R_u\cdots R_u
[\D_u,a_j]R_u\dot{\D_u}R_u[\D_u,a_{j+1}]\cdots R_u[\D_u,a_M]R_u.
\end{eqnarray*}
\end{lemma}

\begin{proof}
First we employ the estimate for $a^*=-a$, which is equation 10.58 in
\cite{GVF},
$$ -\Vert [\D,a]\Vert\dd^{-u}\leq [\D_u,a]\leq \Vert [\D,a]\Vert\dd^{-u} $$
to deduce that $[\D_u,a]\in\LL^{q^+/u}$ where we write $q^+$ to indicate
$q+\epsilon$ for all positive $\epsilon.$ Using $\D_u\in OP^{1-u}$, we
deduce that $R_u\in OP^{-2+2u}$ and so $R_u\in \LL^{{q^+}/2(1-u)}$. The
operator $T_{u,j}$ has $M+2$ factors of $R_u$, one factor of $\D_u$, and
$M+1$ factors $[\D_u,a_j]$, and so $T\in \LL^r$ where
$$\frac{1}{r}=(M+1)\frac{u}{q^+}+(M+3/2)
\frac{2(1-u)}{q^+}=\frac{2M+3-(M+2)u}{q^+}.$$
The worst possible case is when $u=1$ and we find $1/r\geq (M+1)/q^+>1$,
and so $T\in\LL^1$ for all $u\in[0,1]$. The strict inequality follows
since $M=2N-P$, and $q-1<2N-P\leq q+1$. 

The strict inequality also implies the following. There exists
$\rho>0$ such that for all $u\in[0,1]$, $T_{u,j}\dd^\rho\in\LL^1$. We
fix some such choice of $\rho$ from now on. Now, $[\D_u,a_j]\dd^u$
is order zero for all $u$ and uniformly $\leq\n[\D,a_j]\n$, 
and  $R_u\dd^{-2u}$ is of
order $-2$ uniformly in $u$. We want to write $T_{u,j}$ as a product
of terms of order $0$ and $-2$ (and a single term of order $1$). To
do this, we insert in the expression for $T_{u,j}$
a factor of $|\D|^{-2u}$ next to each individual $R_u$ and a factor of
$|\D|^u$ to the right of each commutator $[\D,a]$: we are then forced to 
insert a factor of $|\D|^u$ to the left of the next commutator. In order to 
make the factors near the second commutator have total order $0$ we are
forced into adding a factor of $|\D|^{-u}$ to the right of the second 
commutator which in turn forces us to insert another factor of $|\D|^u$,
etc. We make another adjustment when we get to the term 
$R_u\D_uR_u=R_u\D|\D|^{-u}R_u$, which force further adjustments after that
term. Thus we get:

\bean T_{u,j}&=&[\D_u,a_0]\dd^u\ R_u\dd^{-2u}\ \dd^u[\D_u,a_1]|\D|^u|\D|^{-u}\
R_u\dd^{-2u}\ \dd^{2u}[\D_u,a_2]\dd^{-u}\cdots\nno
&&\cdots\dd^{ju}[\D_u,a_j]\dd^{-(j-1)u}\ R_u\dd^{-2u}\D\ R_u\dd^{-2u}\
\dd^{(j+2)u}[\D_u,a_{j+1}]\ \dd^{-(j+1)u}\cdots\nno
&&\qquad\cdots\dd^{(M+1)u}[\D_u,a_M]\dd^{-Mu}\ R_u\dd^{-2u}\
\dd^{(M+2)u}\nno
&=&[\D_u,a_0]\dd^u\ R_u\dd^{-2u}R_0^{-1}\ R_0\ \dd^u[\D_u,a_1]\
R_u\dd^{-2u}R_0^{-1}\ R_0\ \dd^{2u}[\D_u,a_2]\dd^{-u}\cdots\nno
&&\cdots\dd^{ju}[\D_u,a_j]\dd^{-(j-1)u}\ 
R_u\dd^{-2u}R_0^{-1}R_0\D\ R_u\dd^{-2u}R_0^{-1}R_0\cdot\nno
&&\cdot\dd^{(j+2)u}[\D_u,a_{j+1}]\ \dd^{-(j+1)u}
\cdots\dd^{(M+1)u}[\D_u,a_m]\dd^{-Mu}\ \dd^{-2u}R_uR_0^{-1}\ R_0\
\dd^{(M+2)u}.\eean
In the first equality we have expressed $T_{u,j}$ as a product of
operators of order zero (the operators 
$\dd^{ku}[\D_u,a_k]\dd^{-(k-1)u}$ and $\dd^{(i-1)u}[\D_u,a_i]\dd^{-iu}$ are 
order zero by the pseudo-differential calculus), operators of order $-2$,
($R_u\dd^{-2u}$) and a final term $\dd^{(M+2)u}$ of order $(M+2)u\leq
M+2$. 

In the second equality we observe that (assuming $\lambda = a+ iv$ where
$a=\delta/2$) there is a
uniform estimate which is a consequence of the functional calculus:
$|\dd^{-2u}R_uR_0^{-1}|\leq 2\delta^{-u}\leq 2\delta^{-1}$. Hence,
$\dd^{-2u}R_uR_0^{-1}$ is order zero independent of the values of $u$, and so
$$T_{u,j}=A_0R_0A_1R_0\cdots A_jR_0\D R_0\cdots A_mR_0\dd^{(M+2)u},$$
where the order zero operators $A_j$ are now $u$-dependent.

Obviously
$$T_{u,j}=T_{u,j}\dd^{M+1-\rho}
\dd^{-M-1+\rho}.$$ 
but by our choice of $\rho$, $\dd^{-M-1+\rho}\in\LL^1$, and
$T_{u,j}\dd^{M+1-\rho}$ is uniformly bounded in $\lambda, s, u$ and,
uniformly in these parameters, of  
order $-\rho$ as a pseudo-differential operator. 

We now want to consider the difference quotients, 
\bean&&
\frac{1}{\epsilon}\left(T_{u+\epsilon,j}\dd^{M+1-\rho}-
T_{u,j}\dd^{M+1-\rho}\right).\eean
Observe that  these differences lie in
$OP^{-\rho}$, and so are bounded. In order to 
take the limit as $\epsilon\to 0$, we require some estimates.  If $t$
lies between $x$ and $1$, we have
$$ \int_1^t s^{-1-\epsilon}ds\leq
\int_1^xs^{-1-\epsilon}ds=\frac{1}{\epsilon} (1-x^{-\epsilon}).$$

So let $x\in[\sqrt{\delta},\infty)$, where $\D^2\geq \delta$, and observe that
$\int_1^x t^{-1-\epsilon}dt =-(x^{-\epsilon}-1)/\epsilon.$
Then we have the estimate
\bean
\int_1^x(t^{-1}-t^{-1-\epsilon})dt&=&-\int_1^xt^{-1}(t^{-\epsilon}-1)dt\nno
&=&-\int_1^xt^{-1}(-\epsilon\int_1^ts^{-1-\epsilon}ds)dt\nno
&\leq&\epsilon\int_1^xt^{-1}\int_1^xs^{-1-\epsilon}dsdt\nno
&=&(1-x^{-\epsilon})\log x.\eean
Hence
$$ \left|\frac{1}{\epsilon}(x^{-\epsilon}-1)+\log x\right| \leq
(1-x^{-\epsilon})\log x.$$
Since $|\log x|\leq
C x^\rho$ for all $x\geq\delta$, for any $\rho>0$, the functional
calculus 
then gives us the inequality
$$ \Vert(\frac{1}{\epsilon}(\dd^{-\epsilon}-1)+\log
\dd)\dd^{-\rho}\Vert\leq C\Vert 1-\dd^{-\epsilon}\Vert\Vert
\dd^{-\rho}\log\dd\Vert \to 0\ \mbox{as}\ \epsilon\to 0.$$
Thus 
$$ \lim_{\epsilon\to 0}\frac{1}{\epsilon}(\dd^{-\epsilon}-1)\dd^{-\rho}
=-\dd^{-\rho}\log\dd $$
where the limit is taken in the norm topology. 
Since we are using the norm topology, and both sides lie in
$OP^{0}$, we may regard this as an equality in $OP^{0}$.

We will now show that the limit
$$ \lim_{\epsilon\to
0}\frac{1}{\epsilon}(T_{(u+\epsilon),j}-T_{u,j})\dd^{M+1-\rho}$$
also exists in $OP^0$ in the norm topology, and therefore.
$\lim_{\epsilon\to
0}\frac{1}{\epsilon}(T_{(u+\epsilon),j}-T_{u,j})$
exists in the $\mathcal L^1$ norm.
Our earlier computations show that we can write 
$$ T_{u,j}=B_0(u)\br_uB_1(u)\br_u\cdots B_j(u)\br_u\D\br_u\cdots
\br_uB_M(u)\br_u ,$$
where $\br_u=R_u\dd^{-2u}\in OP^{-2}$ and for each $k=0,\dots,M$
$$ B_k(u)=\dd^{ku}[\D_u,a_k]\dd^{-(k-1)u}\ \ \mbox{or}\ \
B_k(u)=\dd^{(k-1)u}[\D_u,a_k]\dd^{-ku},$$
depending on whether $k\leq j$ or $k>j.$ 
and in each case $B_k(u)\in OP^0$. So now we can write
\bean &&\frac{1}{\epsilon}(T_{u+\epsilon}-T_u)\dd^{M+1-\rho}\nno
&=&\frac{1}{\epsilon}(B_0(u+\epsilon)-B_0(u))\dd^{-\rho} \dd^\rho
\br_{u+\epsilon} B_1(u+\epsilon)\cdots
B_M(u+\epsilon)\br_{u+\epsilon}\dd^{M+1 -\rho}\nno 
&+& \frac{1}{\epsilon}B_0(u)(\br_{u+\epsilon}-\br_u) \dd^{-\rho}
\dd^\rho B_1(u+\epsilon)\br_{u+\epsilon}\cdots
B_M(u+\epsilon)\br_{u+\epsilon}\dd^{M+1 -\rho}\nno 
&+&\dots\nno
&+& \frac{1}{\epsilon} B_0(u)\br_u B_1(u)\cdots
B_M(u)(\br_{u+\epsilon}-\br_u)\dd^{-\rho} \dd^\rho\dd^{M+1 -\rho}. \eean
We have written each term as a product $\alpha\beta\gamma$, where the
sum of orders of $\alpha$ and $\gamma$ is two or zero, while $\beta$ is a
difference quotient times $\dd^{-\rho}$, and so the order of $\beta$
is always either $-2-\rho$ or $-\rho$. Below we show that the two
possiblities for $\beta$ give norm convergent limits, in $OP^{-2}$ or
$OP^0$, 
which will imply
that $T_{u,j}$
is differentiable, by a standard argument. The continuity of the
resulting derivative can be determined by similar, but simpler,
arguments.

So we now examine the difference quotients.
There are two kinds of terms which arise in this computation (the
factor of $\dd^{-u}$ arising from the extra $\D_u$ is included in our
other factors). First we consider $\br_{u+\epsilon}-\br_u$ which equals:
\bean R_{u+\epsilon}\dd^{-2(u+\epsilon)}-R_u\dd^{-2u}&=&
R_{u+\epsilon}\dd^{-2(u+\epsilon)}-R_u\dd^{-2(u+\epsilon)}+
R_u\dd^{-2(u+\epsilon)}-R_u\dd^{-2u} \nno
&=&R_{u+\epsilon}\dd^{-2(u+\epsilon)}(\dd^{-2\epsilon}-1)\D_u^2R_u +
R_u\dd^{-2u}(\dd^{-2\epsilon}-1).\eean
Observe we have chosen to add and subtract a term which lies in
$OP^{-2}$. 
Since the factor $R_u\dd^{-2u}=R_u\dd^{-2u}R_0^{-1}R_0$ 
is uniformly in $OP^{-2}$, the same method of proof as before shows
that
$$\lim_{\epsilon\to 0} \frac{1}{\epsilon}(
R_{u+\epsilon}\dd^{-2(u+\epsilon)}-R_u\dd^{-2u})\dd^{-\rho}
=R_u\dd^{-2u}2\log\dd(\D_u^2R_u-1)\dd^{-\rho} $$
where the limit exists in norm, and the limit lies in
$OP^{-2}$.

The second kind of difference quotient is (one of the two forms of) 
$B_n(u+\epsilon)-B_n(u)$, and we will present the proof for the case 
$$ B_n(u)=\dd^{nu}[\D_u,a_n]\dd^{-(n-1)u},$$
the other possibility being entirely similar. By adding and subtracting terms 
of order strictly less than zero, we can rewrite this difference as

\bean &&\dd^{n(u+\epsilon)}[\D_{u+\epsilon},a]\dd^{-(n-1)(u+\epsilon)}-
\dd^{nu}[\D_u,a]\dd^{-(n-1)u}\nno
&=&
\dd^{n(u+\epsilon)}[\D_{u+\epsilon},a]\dd^{-(n-1)(u+\epsilon)}-
\dd^{n(u+\epsilon)}[\D_{u+\epsilon},a]\dd^{-(n-1)u}\nno
&&+
\dd^{n(u+\epsilon)}[\D_{u+\epsilon},a]\dd^{-(n-1)u}-\dd^{nu}[\D_{u+\epsilon},a]
\dd^{-(n-1)u} \nno
&& +\dd^{nu}[\D_{u+\epsilon},a]\dd^{-(n-1)u}-\dd^{nu}[\D_u,a]\dd^{-(n-1)u}\nno
&=&\dd^{n(u+\epsilon)}[\D_{u+\epsilon},a]\dd^{-(n-1)u}(\dd^{-(n-1)\epsilon}-1) + 
(1-\dd^{-n\epsilon})\dd^{n(u+\epsilon)}[\D_{u+\epsilon},a]\dd^{-(n-1)u}\nno
&&+\dd^{nu}[\D_{u+\epsilon}-\D_u,a]\dd^{-(n-1)u}.\eean
The first term can be handled using our estimates above, the second term can
be handled similarly by using the factor of $|\D|^{\rho} |\D|^{-\rho}$ 
on the left
so only the last term needs examination. We rewrite this last term as 
\bean &&\dd^{nu}[\D_{u+\epsilon}-\D_u,a]\dd^{-(n-1)u}\nno
&=&\dd^{nu}\D[\dd^{-(u+\epsilon)}-\dd^{-u},a]\dd^{-(n-1)u}+\dd^{nu}[\D,a]
\dd^{-nu}(\dd^{-\epsilon}-1).\eean
The second of these terms can also be dealt with using our previous methods, 
and so we are left with examining
$$\frac{1}{\epsilon}\dd^{nu}\D[(\dd^{-(u+\epsilon)}-\dd^{-u}),a]\dd^{-(n-1)u}.$$
We proceed in stages. 
Recalling that this difference 
quotient is multiplied 
by $\dd^{-\rho}$, we  consider the convergence of
$$ \frac{1}{\epsilon}\dd^{nu}\D[(\dd^{-(u+\epsilon)}-\dd^{-u}),a]
\dd^{-\rho}\dd^{-(n-1)u}$$
in $OP^{-\rho}\subset OP^0$. 

Applying the Leibnitz rule we have 
\begin{align*}
\frac{1}{\epsilon}\dd^{nu}\D[(\dd^{-(u+\epsilon)}-\dd^{-u}),a]
\dd^{-\rho}\dd^{-(n-1)u}&= 
\frac{1}{\epsilon}\dd^{(n-1)u}\D[(\dd^{-\epsilon}-1),a]\dd^{-\rho}\dd^{-(n-1)u}\nno
&+ 
\frac{1}{\epsilon}\dd^{nu}\D[\dd^{-u},a](\dd^{-\epsilon}-1)
\dd^{-\rho}\dd^{-(n-1)u}.\end{align*}
As $\dd^u\D[\dd^{-u},a]\in OP^0$, and conjugation by $\dd^{(n-1)u}$
preserves $OP^0$, we see using our previous methods that the second
term has a limit in $OP^0$.
For  
the first term we employ the Leibnitz rule again. This gives us
\begin{align*}&\frac{1}{\epsilon}\dd^{(n-1)u}\D[(\dd^{-\epsilon}-1),a]
\dd^{-\rho}\dd^{-(n-1)u}\nno
& = 
\frac{1}{\epsilon}\dd^{(n-1)u}\D[(\dd^{-\epsilon}-1)\dd^{-\rho},a]
\dd^{-(n-1)u}\nno
&+ 
\frac{1}{\epsilon}\dd^{(n-1)u}(\dd^{-\epsilon}-1)\dd^{-\rho}\D[\dd^{\rho},a] 
\dd^{-\rho}\dd^{-(n-1)u}.\end{align*}
As before conjugation by $\dd^{(n-1)u}$ does not affect matters. For
the first term we observe that
$\D[\frac{1}{\epsilon}(\dd^{-\epsilon}-1)\dd^{-\rho},a]$ is uniformly
in $OP^0$ and has a limit as $\epsilon\to 0$ by our previous
methods. The proof is completed
by noting that the second term is handled similarly once we  see that 
$$\D[\dd^{\rho},a]\dd^{-\rho}=-\D\dd^{\rho}[\dd^{-\rho},a]\in OP^0.$$
\end{proof}

\begin{proof}(of Lemma \ref{aagghhh}).
Lemma \ref{had-to-be-done}, and together with arguments
of a similar nature, show that $\Psi^r_{u,M}$
and $\frac{d}{du}\Phi^r_{u,M+1}$ are well-defined and are
continuous. The proof of Lemma \ref{had-to-be-done} also shows that 
the formal differentiations given below are in fact justified.

First of all, using the $\D_u$ version of
Equation \ref{formulaforBPhi} of Lemma \ref{tcohom} and the $R_u$ version of 
Definition \ref{newexpectation} to expand $(B\Phi^r_{u,M+1})(a_0,...,a_M)$,
we see that it is the sum of the $T_{u,j}$ and so its derivative is the
sum over $j$ of the derivatives in Lemma \ref{had-to-be-done}. Using the $R_u$
version of
Definition \ref{newexpectation} again to rewrite this in terms of 
$\la\la\cdots\ra\ra$ where possible we get:
\bean &&\frac{d}{du}(B\Phi^r_{u,M+1})(a_0,...,a_M)\nno
&=&-\frac{\eta_M}{2}\int_0^\infty s^M\sum_{i=0}^M\left(\la\la 
[\D_u,a_0],...,[\D_u,a_i],2\D_u\dot{\D}_u,...,[\D_u,a_M]\ra\ra_{M+1,s,r}\right.
\nno
&&+\left.\la\la 
[\D_u,a_0],...,[\dot{\D}_u,a_i],...,[\D_u,a_M]\ra\ra_{M,s,r}\right) ds\nno
&&-\frac{\eta_M}{2}\int_0^\infty s^M\sum_{i=0}^M(-1)^i\la 
[\D_u,a_0],...,[\D_u,a_i],\dot{\D}_u,...,[\D_u,a_M]\ra_{M+1,s,r} ds.\eean
For the next step we compute $Bb\Psi^r_{u,M}$, and then use $bB=-Bb$. First we
apply $b$
\bean &&(b\Psi^r_{u,M})(a_0,...,a_{M+1})\nno
&=&-\frac{\eta_M}{2}\sum_{j=1}^M(-1)^j\int_0^\infty s^M\la\la 
a_0\dot{\D}_u,...,[\D_u,a_ja_{j+1}],...,[\D_u,a_{M+1}]\ra\ra_{M,s,r} ds\nno
&&-\frac{\eta_M}{2}\int_0^\infty s^M\la\la 
a_0a_1\dot{\D}_u,[\D_u,a_2],...,[\D_u,a_{M+1}]\ra\ra_{M,s,r} ds\nno
&&-(-1)^{M+1}\frac{\eta_M}{2}\int_0^\infty s^M\la\la 
a_{M+1}a_0\dot{\D}_u,[\D_u,a_1],...,[\D_u,a_M]\ra\ra_{M,s,r} ds\eean
\bean
&=&-\frac{\eta_M}{2}\sum_{j=1}^M(-1)^j\int_0^\infty s^M\la\la 
a_0\dot{\D}_u,...,a_j[\D_u,a_{j+1}]+[\D_u,a_j]a_{j+1},...,[\D_u,a_{M+1}]\ra\ra_{
M,s,r} ds\nno
&&-\frac{\eta_M}{2}\int_0^\infty s^M\la\la 
a_0\dot{\D}_ua_1,[\D_u,a_2],...,[\D_u,a_{M+1}]\ra\ra_{M,s,r} ds\nno
&&+\frac{\eta_M}{2}\int_0^\infty s^M\la\la 
a_0[\dot{\D}_u,a_1],[\D_u,a_2],...,[\D_u,a_{M+1}]\ra\ra_{M,s,r} ds\nno
&&-(-1)^{M+1}\frac{\eta_M}{2}\int_0^\infty s^M\la\la 
a_{M+1}a_0\dot{\D}_u,[\D_u,a_1],...,[\D_u,a_M]\ra\ra_{M,s,r} ds\eean
\bean
&=&-\frac{\eta_M}{2}\int_0^\infty s^M\sum_{j=1}^{M+1}(-1)^j\la\la 
a_0\dot{\D}_u,[\D_u,a_1],...,[\D^2_u,a_j],...,[\D_u,a_{M+1}]\ra\ra_{M+1,s,r} 
ds\nno
&&-\frac{\eta_M}{2}\int_0^\infty 
s^M\sum_{j=1}^{M+1}(-1)^j(-1)^{deg(a_0\dot{\D}_u)+\cdots+deg([\D_u,a_{j-1}])}\la 
a_0\dot{\D}_u,[\D_u,a_1],...,[\D_u,a_{M+1}]\ra_{M+1,s,r} ds\nno
&&+\frac{\eta_M}{2}\int_0^\infty s^M\la\la 
a_0[\dot{\D}_u,a_1],...,[\D_u,a_{M+1}]\ra\ra_{M,s,r} ds\eean
The last equality follows from the 
$R_u$ version of Lemma \ref{lala-la-identity}. In the above we note that
$deg(a_0\dot{\D_u})=1= deg([\D_u,a_k])$ for all $k$ so that 
$deg(a_0\dot{\D_u})+\cdots+deg([\D_u,a_{j-1}])=j$ and 
$deg(a_0\dot{\D_u})+\cdots+deg([\D_u,a_{M+1}])=M+2\equiv P(mod\ 2).$
We also note the commutator identity $[\D_u^2,a_j]=\{\D_u,[\D_u,a_j]\}=[\D_u,[\D_u,a_j]]_{\pm}$
so in order to apply the $\D_u$ version of Equation (6) of Lemma 
\ref{lala-la-identity} we first add and subtract:
\ben -\frac{\eta_M}{2}\int_0^\infty s^M \la\la 
\{\D_u,a_0\dot{\D}_u\},[\D_u,a_1],...,[\D_u,a_{M+1}]\ra\ra_{M+1,s,r} ds\een
and then apply Equation (6) to get:
\bean &&-2\frac{\eta_M}{2}\int_0^\infty s^M \sum_{j=0}^{M+1}\la 
a_0\dot{\D}_u,...,[\D_u,a_j],\D^2_u,...,[\D_u,a_{M+1}]\ra_{M+2,s,r} ds\nno
&&+\frac{\eta_M}{2}\int_0^\infty s^M\la\la 
a_0\{\D_u,\dot{\D}_u\}+[\D_u,a_0]\dot{\D}_u,[\D_u,a_1],...,[\D_u,a_{M+1}]\ra\ra_
{M+1,s,r} ds\nno
&&-\frac{\eta_M}{2}(M+1)\int_0^\infty s^M \la 
a_0\dot{\D}_u,[\D_u,a_1],...,[\D_u,a_{M+1}]\ra_{M+1,s,r} ds\nno
&&+\frac{\eta_M}{2}\int_0^\infty s^M \la\la 
a_0[\dot{\D}_u,a_1],...,[\D_u,a_{M+1}]\ra\ra_{M,s,r} ds\eean
Then we apply the $\D_u$ version of Lemma \ref{differentfort} as modified in the
proof of Proposition \ref{tcohom} with $t=0$ to the first term above to get
\bean
&=&\frac{\eta_M}{2}(q +2r)\int_0^\infty s^M\la 
a_0\dot{\D}_u,[\D_u,a_1],...,[\D_u,a_{M+1}]\ra_{M+1,s,r} ds\nno
&&+\frac{\eta_M}{2}\int_0^\infty s^M\la\la 
a_0\{\D_u,\dot{\D}_u\}+[\D_u,a_0]\dot{\D}_u,[\D_u,a_1],...,[\D_u,a_{M+1}]\ra\ra_
{M+1,s,r} ds\nno
&&+\frac{\eta_M}{2}\int_0^\infty s^M \la\la 
a_0[\dot{\D}_u,a_1],...,[\D_u,a_{M+1}]\ra\ra_{M,s,r} ds\eean
The next step is to apply $B$ to these three terms:
\bean &&(Bb\Psi^r_{u,M})(a_0,...,a_M)\nno
&=&(q+2r)\frac{\eta_M}{2}\sum_{j=0}^M(-1)^{(M+1)j}\int_0^\infty s^M\la 
\dot{\D}_u,[\D_u,a_j],...,[\D_u,a_{j-1}]\ra_{M+1,s,r} ds\nno
&&+\frac{\eta_M}{2}\sum_{j=0}^M(-1)^{(M+1)j}\int_0^\infty s^M \la\la 
\{\D_u,\dot{\D}_u\},[\D_u,a_j],...,[\D_u,a_{j-1}]\ra\ra_{M+1,s,r} ds\nno
&&+\frac{\eta_M}{2}\sum_{j=0}^M(-1)^{(M+1)j}\int_0^\infty s^M\la\la 
[\dot{\D}_u,a_j],...,[\D_u,a_{j-1}]\ra\ra_{M,s,r} ds\eean
\bean
&=&\frac{(q+2r)\eta_M}{2}\sum_{j=0}^M(-1)^{(M+1)j+Aj}\int_0^\infty s^M\la 
[\D_u,a_0],...,[\D_u,a_{j-1}],\dot{\D}_u,...,[\D_u,a_M]\ra_{M+1,s,r} 
ds\nno
&+&\frac{\eta_M}{2}\sum_{j=0}^M(-1)^{(M+1)j+(1+A)j}\int_0^\infty s^M 
\la\la[\D_u,a_0],..., 
\{\D_u,\dot{\D}_u\},[\D_u,a_j],...,[\D_u,a_M]\ra\ra_{M+1,s,r} ds\nno
&+&\frac{\eta_M}{2}\sum_{j=0}^M(-1)^{(M+1)j+(1+A)j}\int_0^\infty 
s^M\la\la[\D_u,a_0],...,[\D_u,a_{j-1}], 
[\dot{\D}_u,a_j],...,[\D_u,a_M]\ra\ra_{M,s,r} ds\eean
\bean
&=&(q+2r)\frac{\eta_M}{2}\sum_{j=0}^M(-1)^j\int_0^\infty s^M\la 
[\D_u,a_0],...,[\D_u,a_{j-1}],\dot{\D}_u,...,[\D_u,a_M]\ra_{M+1,s,r} 
ds\nno
&&+\frac{\eta_M}{2}\sum_{j=0}^M\int_0^\infty s^M \la\la[\D_u,a_0],..., 
2\D_u\dot{\D}_u,[\D_u,a_j],...,[\D_u,a_M]\ra\ra_{M+1,s,r} ds\nno
&&+\frac{\eta_M}{2}\sum_{j=0}^M\int_0^\infty 
s^M\la\la[\D_u,a_0],...,[\D_u,a_{j-1}], 
[\dot{\D}_u,a_j],...,[\D_u,a_M]\ra\ra_{M,s,r} ds\eean
Using $bB=-Bb$, and our formula for 
$\frac{d}{du}(B\Phi_{u,M+1}^r)(a_0,...,a_M)$ we get:
\bean &&(bB\Psi^r_{u,M})(a_0,...,a_M)\nno
&=&-(q+2r)\frac{\eta_M}{2}\sum_{j=0}^M(-1)^j\int_0^\infty s^M\la 
[\D_u,a_0],...,[\D_u,a_{j-1}],\dot{\D}_u,...,[\D_u,a_M]\ra_{M+1,s,r} 
ds\nno
&&+\frac{\eta_M}{2}\sum_{i=0}^M(-1)^i\int_0^\infty s^M \la 
[\D_u,a_0],...,[\D_u,a_i],\dot{\D}_u,...,[\D_u,a_M]\ra_{M+1,s,r} ds\nno
&&+\frac{d}{du}(B\Phi_{u,M+1}^r)(a_0,...,a_M).\eean
This proves the result.
\end{proof}

\begin{cor}\label{approxderiv} For $\D$ invertible we have
\be 
\frac{1}{(r+(q-1)/2)}(bB\Psi^r_{u,M})(a_0,...,a_M)=\frac{1}{(r+(q-1)/2)}\frac{d}
{du}(B\Phi^{r}_{u,M+1})(a_0,...,a_M)+h(r)\label{deriveqn}\ee
where $h(r)$ is analytic for $Re(r)>-M/2$.
\end{cor}

\begin{proof} This follows from \cite[Lemma 7.2]{CPRS2} applied to the function
\ben \int_0^\infty s^M\la 
[\D_u,a_0],...,[\D_u,a_i],\dot{\D}_u,...,[\D_u,a_M]\ra_{M+1,s,r} ds,\een
which shows that this function is holomorphic for $Re(r)>-M/2$, and in both the 
even and odd cases, $-M/2<(1-q)/2$.
\end{proof}

\begin{cor}\label{exactderiv} For $\D$ invertible we have
\ben 
(bB\Psi^{(1-q)/2}_{u,M})(a_0,...,a_M)=\frac{d}{du}(B\Phi^{(1-q)/2}_{u,M+1})(a_0,
...,a_M).\een
\end{cor}

\begin{proof} We have already observed that $B\Phi^r_{u,M+1}$ is holomorphic at 
$r=(1-q)/2$, and so by Corollary \ref{approxderiv}, we can take the residues of 
both sides of Equation (\ref{deriveqn}) to obtain the result. Observe that 
taking these residues did not require isolated spectral dimension.
\end{proof} 

\subsection{Proofs of the Theorems \ref{approx} and \ref{exactly}.}
 
\begin{proof}({\em Theorem \ref{approx}})
The image of the cyclic cochain
\ben \frac{1}{(r-(1-q)/2)}\int_0^1B\Psi^{r}_{u,M}(a_0,...,a_M)du\een
under the operator  $b$ is given by
\be 
\frac{1}{(r-(1-q)/2)}\int_0^1\frac{d}{du}B\Phi^{r}_{u,M+1}(a_0,...,a_M)du+ho
lo.\label{intofderiv}\ee
Here $holo$ is the integral over $u$ of the holomorphic remainder from Corollary 
\ref{approxderiv}. Integrating this remainder in $u$ does not affect the 
estimates proving holomorphicity at $r=(1-q)/2$, since the integral is 
absolutely convergent. By the fundamental theorem of calculus, Equation 
(\ref{intofderiv}) is (modulo functions holomorphic at $r=(1-q)/2$) the 
difference of $\frac{1}{(r-(1-q)/2)}B\Phi^{r}_{M+1}$ defined using 
$F=\D\dd^{-1}$ and $\frac{1}{(r-(1-q)/2)}B\Phi^{r}_{M+1}$ defined using $\D$. 
Hence the two are cohomologous in cyclic cohomology.  Recalling that $F^2=1$ and 
using our previous formula for $B\Phi_{u,m}^r$ (the $\D_u$ version of
Proposition \ref{tcohom} with $u=1$) we have
\bean && (B\Phi^{r}_{M+1})(a_0,...,a_M)|_{u=1}\nno
&=&-\frac{\eta_{M}}{2}\sum_{j=0}^{M}(-1)^{j+1}\int_0^\infty s^{M}\la 
[F,a_0],...,[F,a_j],F,[F,a_{j+1}],...,[F,a_M]\ra_{M+1,s,r} ds\nno
&=&-\frac{\eta_{M}}{2}\sum_{j=0}^M\int_0^\infty s^{M}\frac{1}{2\pi 
i}\tau\left(\gamma\int_l\lambda^{-q/2-r}F[F,a_0]\cdots[F,a_M](\lambda-(s^2+1))^
{-M-2}d\lambda\right)ds\nno
&=&\frac{\eta_{M}}{2}\frac{(-1)^{M}}{M!}\frac{\Gamma(M+1+q/2+r)}{\Gamma(q/2+r)}\
\int_0^\infty s^{M}\tau(\gamma F[F,a_0]\cdots[F,a_M](s^2+1)^{-M-1-q/2-r})ds\eean
In the second equality we anticommuted $F$ past the commutators, and pulled all 
the resolvents to the right (they commute with everything, since they involve 
only scalars.) In the last equality we used the Cauchy integral formula to do 
the contour integral, and performed the sum. 

Now we pull out $(s^2+1)^{-M-1-p/2-r}$ from the trace, leaving the identity 
behind. The $s$ integral is as follows.
\bean &&\int_0^\infty s^M(s^2+1)^{-M-1-p/2-r}ds\nno
&=&\frac{1}{\Gamma(M+1+q/2+r)}\int_0^\infty\int_0^\infty 
s^Mu^{M+q/2+r}e^{-u(s^2+1)}duds\nno
&=&\frac{1}{\Gamma(M+1+q/2+r)}\int_0^\infty\int_0^\infty 
s^Mu^{M+q/2+r}e^{-u(s^2+1)}dsdu\nno
&=&\frac{\Gamma((M+1)/2)}{2\Gamma(M+1+q/2+r)}\int_0^\infty 
u^{q/2+r+M/2-1/2}e^{-u}du\nno
&=&\frac{\Gamma((M+1)/2)\Gamma(q/2+r+M/2+1/2)}{2\Gamma(M+1+1/2)}\eean
Putting the pieces together gives
\bean &&(B\Phi^{r}_{M+1})(a_0,...,a_M)|_{u=1}\nno
&=&\frac{\eta_M}{2}(-1)^M\frac{\Gamma((M+1)/2)}{\Gamma(q/2+r)}\frac{\Gamma(((q-1
)/2+r)+M/2+1)}{2M!}\tau(\gamma F[F,a_0]\cdots[F,a_M])\eean
Now  $\eta_M=\sqrt{2i}^P(-1)^M2^{M+1}\Gamma(M/2+1)/\Gamma(M+1)$, and the 
duplication formula for the Gamma function tells us that
\ben \Gamma((M+1)/2)\Gamma(M/2+1)2^M=\sqrt{\pi}\Gamma(M+1).\een
Hence 
\ben 
(B\Phi^{r}_{M+1})(a_0,...,a_M)|_{u=1}=\frac{\sqrt{\pi}\sqrt{2i}^P
\Gamma(((q-1)/2
+r)+M/2+1)}{\Gamma(q/2+r)2\cdot M!}\tau(\gamma 
F[F,a_0][F,a_1]\cdots[F,a_M]).\een
Now we use the functional equation for the Gamma function
\bean &&\Gamma(((q-1)/2+r)+M/2+1)\nno
&&=\Gamma((q-1)/2+r)\times((q-1)/2+r+M/2)\times((q-1)/2+r+M/2-1)
\cdots((q-1)/2+r
)\eean
to write this as
\ben (B\Phi^{r}_{M+1})(a_0,...,a_M)|_{u=1}=\frac{C_{q/2+r}\sqrt{2i}^P}{2\cdot 
M!}\sum_{j=1}^{M/2}((r+(q-1)/2)^j\s_{M/2,j}\tau(\gamma 
F[F,a_0][F,a_1]\cdots[F,a_M]),\een
where the $\s_{M/2,j}$ are elementary symmetric functions of the integers 
$1,2,...,M/2$ (even case) or of the half integers $1/2,3/2,\dots,M/2$ 
(odd case). Recalling that the `constant' $C_{q/2+r}$ has a simple pole at 
$r=(1-q)/2$ with residue equal to $1$, and $\s_{M/2,1}=\Gamma(M/2+1)$ in 
both even and odd cases, and recalling Definition \ref{conditional} of
$\tau^\prime$ we see that
\bean &&\frac{1}{(r-(1-q)/2)}(B\Phi^{r}_{M+1})(a_0,...,a_M)|_{u=1}\nno
&=&\frac{\sqrt{2i}^P\Gamma(M/2+1)}{(r-(1-q)/2)2\cdot M!}\tau(\gamma 
F[F,a_0]\cdots[F,a_M])+holo\nno
&=&\frac{\sqrt{2i}^P\Gamma(M/2+1)}{M!(r-(1-q)/2)}\tau'(\gamma 
a_0[F,a_1]\cdots[F,a_M])+holo\nno
&=&\frac{1}{(r-(1-q)/2)}Ch_F(a_0,a_1,\dots,a_M)+holo,\eean
where $holo$ is a function holomorphic at $r=(1-q)/2$, and on the right 
hand side the Chern character appears with its $(b,B)$ normalisation.
\end{proof}

\begin{proof}({\em Theorem \ref{exactly}}). If we assume isolated spectral 
dimension we can take residues of the resolvent cocycle to obtain the residue 
cocycle. By Corollary \ref{isolatedcohom}, the residue cocycle is cohomologous 
to $B\Phi^{(1-q)/2}_{M+1}$. Observe that it is only in Corollary 
\ref{isolatedcohom} that we need to assume isolated spectral dimension. This is 
because by Corollary \ref{exactderiv} we always have 
\ben 
(bB\Psi^{(1-q)/2}_{u,M})(a_0,...,a_M)=\frac{d}{du}(B\Phi^{(1-q)/2}_{u,M+1})(a_0,
...,a_M).\een
Then by the computations above in the proof of Theorem \ref{approx}, we have 
\ben (B\Phi^{(1-q)/2}_{u=1,M+1})(a_0,...,a_M)\mbox{ is cohomologous to 
}(B\Phi^{(1-q)/2}_{u=0,M+1})(a_0,...,a_M)\een
which again by the proof of Theorem \ref{approx} completes the proof of the 
Theorem.
\end{proof}

\subsection{Removing the Invertibility Hypothesis}\label{noninvtble}

Theorem \ref{res=chern} will have been proved once we show how to remove the 
invertibility
hypothesis.
We shall employ the `double' of the spectral triple $(\A,\HH,\D)$ from Section 
\ref{cycliccohomology}.

In the double up procedure we will start with $0\leq \mu<1$.
We are interested in the relationship between $1+\D^2$ (implicitly
tensored by $Id_2$ here and below) and $1+\D_\mu^2$, given 
by
\ben 1+\D_\mu^2=\bma 1+\mu^2+\D^2 & 0\\ 0 & 1+\mu^2+\D^2\ema.\een
If we scale  $\D_\mu$ by $(1-\mu^2)^{-1/2}$ then we get 
\ben (1+\D_\mu^2)^{-s}\tto (1-\mu^2)^{s}(1+\D^2)^{-s}.\een

Let $\omega_{m,k}=a_0[\D,a_1]^{(k_1)}\cdots[\D,a_m]^{(k_m)}$. Then if we scale 
$\D$ by $\epsilon$, $\omega_{m,k}\to\epsilon^{2|k|+m}\omega_{m,k}$.
If we write $\omega_{\mu,m,k}$ for $\omega_{m,k}$ defined using $\D_\mu$, then 
\ben \omega_{\mu,m,k}=\omega_{m,k}+O(\mu),\een
where the $O(\mu)$ term is an operator of order $|k|$. In terms of the matrix 
representation we have
\ben \omega_{\mu,m,k}=\bma \omega_{m,k} + \mu\omega'_{m,k} & \mu\omega''_{m,k}\\ 
0 & 0\ema.\een

Now take $\omega_{\mu,m,k}(1+\D_\mu^2)^{-r/2-|k|-m/2+1/2-p/2}$ and scale 
$\D_\mu$ by $(1-\mu^2)^{-1/2}$. We obtain
\ben 
(1-\mu^2)^{(q-1)/2+r}\omega_{m,k}(1+\D^2)^{-r-|k|-m/2-(q-1)/2}+O(\mu)
(1+\D^2)^{-
r-|k|-m/2-(q-1)/2},\een
where again the $O(\mu)$ term is an order $|k|$ operator. 

Let us write $\zeta^r_{m,\mu}$ for the sum of zeta functions 
we get by performing the pseudodifferential calculus on the resolvent cocycle 
and discarding the holomorphic remainder. Then, modulo coboundaries and 
functions holomorphic at 
$r=(1-q)/2$,
\bean \frac{1}{(r-(1-q)/2)}Ch_{F_\mu}&=&(\phi^r_{m,\mu})_{m=P}^M\nno
&=&(\zeta^r_{m,\mu})_{m=P}^M\nno
&=&(1-\mu^2)^{r-(1-q)/2}(\zeta^r_{m,0})_{m=P}^M+O(\mu)\nno
&=&(1-\mu^2)^{r-(1-q)/2}(\phi^r_m)_{m=P}^M+O(\mu)\nno
&=&(\phi^r_m)_{m=P}^M+ O(\mu)\eean
where each $O(\mu)$ is a cocycle  with the same regularity properties as the 
resolvent cocycle, and which is zero at $\mu=0$. 

Evaluating both sides of this equation on a $b,B$-cycle yields  functions 
defined and holomorphic in some half plane, and the two sides differ by 
functions holomorphic in a half-plane containing $r=(1-q)/2$. The left hand 
side yields a function  independent of $\mu$ 
by Proposition \ref{chind}, and so either the $O(\mu)$ contributions are 
coboundaries (and so vanish when evaluated on a cycle), or they are 
holomorphic at $r=(1-q)/2$. In either case we find that modulo cochains 
yielding functions holomorphic at $r=(1-q)/2$,
$$ (\phi^r_m)_{m=P}^M \ \ \mbox{is cohomologous to}\ \ \frac{1}{(r-(1-q)/2)}
Ch_{F_\mu}.$$
Taking residues, in the case that $(\A,\HH,\D)$ has isolated spectral dimension,
leads to the analogous result for the residue cocycle:
$$(\phi_m)_{m=P}^M\ \ \ \mbox{is cohomologous to}\ \ Ch_{F_\mu}.$$


Since  $Ch(\A,\HH,\D)$ is $Ch_{F_\mu}$ which is $Ch(\A,\HH^2,\D_{\mu})$ for 
any positive $\mu$, we are done.


\begin{thm} If $(\A,\HH,\D)$ is a $QC^\infty$ finitely summable spectral triple 
with spectral dimension $q\geq 1$, then the resolvent cocycle is cohomologous to
\ben \frac{1}{(r-(1-q)/2)}Ch_{F_\mu}\een
modulo cochains with values in the functions holomorphic at $r=(1-q)/2$. If 
$(\A,\HH,\D)$ also has isolated spectral dimension, then the cyclic cohomology 
class of the residue cocycle coincides with the class of the Chern character of 
$(\HH,F=\D(1+\D^2)^{-1/2})$.
\end{thm}

\section{Some corollaries and the connection with Higson's cocycle}

The transgression cocycle of the previous Section 
allows us to prove a couple of interesting 
corollaries.

\begin{cor} For any $(b,B)$ cycle $(c_m)_{m=P}^K$ of the same parity as 
$(\A,\HH,\D)$, the function 
\ben \sum_{m=P}^M\phi^r_m(c_m)\een
has an analytic continuation to a deleted neighbourhood of $r=(1-q)/2$ with at 
worst a simple pole at $r=(1-q)/2$.
\end{cor}

\begin{proof} The cocycle $(\phi^r_m)$ differs from $(1/(r-(1-q)/2))Ch_F$ by 
coboundaries and functions holomorphic at $r=(1-q)/2$.
\end{proof}

\begin{cor} For any Hochschild $M$-cycle $c_M$ the function
\ben \phi^r_M(c_M)\een
has an analytic continuation to a deleted neighbourhood of $r=(1-q)/2$ with at 
worst a simple pole at $r=(1-q)/2$. If $M$ is even (resp. odd) and $[q]$ is odd 
(resp. even) the residue at $r=(1-q)/2$ vanishes.
\end{cor}

\begin{proof} We have the formula
\ben b\Phi^r_{M-1}+B\Phi^r_{M+1}=(r-(1-q)/2)\phi^r_M.\een
If $c_M=\sum_ia^i_0\otimes a^i_1\otimes\cdots\otimes a^i_M$ is a Hochschild 
cycle, we have $b\Phi^r_{M-1}(c_M)=0$, so
\ben \frac{1}{r-(1-q)/2}B\Phi^r_{M+1}(c_M)=\phi^r_M(c_M).\een
Since there exists $\delta>0$ such that $B\Phi^r_{M+1}(c_M)$ is holomorphic for 
$Re(r)>(1-q)/2-\delta$, we see that $\phi^r_M(c_M)$ meromorphically continues to 
this region with only a simple pole at $r=(1-q)/2$. The region $Re(r)>(1-M)/2$ 
where $\phi^r_M(a_0,...,a_M)$ is holomorphic is $Re(r)>-[q/2]$. For $[q]=2n$, 
$q=2n+\kappa$, $0\leq \kappa<1$ and
\ben \frac{1-q}{2}=-n+\frac{1-\kappa}{2}>-[q/2]=-n.\een
Similar comments apply when $M$ is even and $[q]$ is odd.
\end{proof}

Higson has a cocycle which is evidently similar to our resolvent cocycle.
An essential difference is that from Higson's cocycle, one derives the
unrenormalised local index theorem. We show here that our resolvent
cocycle naturally gives rise to a `renormalised' version of Higson's
cocycle.

We take our resolvent cocycle, perform the pseudodifferential expansion,
the Cauchy integral and the $s$-integral. This gives (modulo functions
holomorphic at $r=(1-q)/2$)

\bean
&&\phi^r_m(a_0,...,a_m)=
\sum_{|k|=0}^{2N-m-P}C(k)(-1)^{m+|k|}\sqrt{\pi}(-1)^P\sqrt{2i}^P
\frac{\Gamma(|k|+(m-1)/2+q/2+r)}{\Gamma(1+|k|+m)\Gamma(q/2+r)}\times\nno
&&\tau(\gamma
a_0[\D,a_1]^{(k_1)}\cdots[\D,a_m]^{(k_m)}(1+\D^2)^{-q/2-r-|k|-(m-1)/2})\eean

We then put back the Cauchy integral using
\bean &&\tau(\gamma
a_0[\D,a_1]^{(k_1)}\cdots[\D,a_m]^{(k_m)}(1+\D^2)^{-q/2-r-|k|-(m-1)/2})\nno
&&(-1)^{|k|+m}\frac{\Gamma(1+|k|+m)\Gamma(q/2+r-(m+1)/2)}{\Gamma(q/2+r+|k|+(m-1)
/2)}
\times\nno
&&\tau\left(\frac{1}{2\pi
i}\int_l\lambda^{-q/2-r+(m+1)/2}a_0[\D,a_1]^{(k_1)}\cdots[\D,a_m]^{(k_m)}
(\lambda-(1+\D^2))^{-|k|-m-1}d\lambda\right)\eean
and undo the pseudodifferential expansion. By our previous estimates,
these operations affect our function-valued cocycle only by functions
holomorphic at the critical point $r=(1-q)/2$. We obtain the following
equality modulo functions holomorphic at $r=(1-q)/2$:
\bean
&&\phi^r_m(a_0,...,a_m)=(-1)^P\sqrt{2i}^P\frac{\sqrt{\pi}\Gamma(q/2+r-(m+1)/2)}
{\Gamma(q/2+r)}\times\nno
&&\tau\left(\frac{1}{2\pi
i}\int_l\lambda^{-q/2-r+(m+1)/2}a_0R_0(\lambda)[\D,a_1]R_0(\lambda)\cdots
[\D,a_m]R_0(\lambda)d\lambda\right).\eean
We call this new cocycle the reduced resolvent cocycle, and denote its
components by $\psi^r_m$ so that the above equality becomes
\ben \phi^r_m(a_0,...,a_m)=\psi^r_m(a_0,...,a_m)\een
modulo functions holomorphic at $r=(1-q)/2$. The integral defining
$\psi^r_m$ exists for $Re(r)>(1-m)/2$ by our previous estimates. The
argument of the coefficent
\ben \Gamma(q/2+r-(m+1)/2)\een
has positive real part when $Re(r)>m/2+(1-q)/2$, and can be
meromorphically continued.

To compare the reduced resolvent cocycle with Higson's improper cocycle,
we write $z=r-(1-q)/2$. Then, writing $\eta^z_m$ for the components of
Higson's improper cocycle we have
\ben \psi^r_m(a_0,...,a_m)=\frac{\sqrt{\pi}}{\Gamma(z+1/2)}\eta^z_m.\een
This gives a `renormalised' version of Higson's cocycle in the sense that
starting with the reduced resolvent cocycle, one arrives at the
renormalised local index theorem, whereas Higson's original cocycle leads
to the unrenormalised theorem.

\end{document}